\documentclass[leqno,11pt]{article}
\usepackage{latexsym}
\usepackage{amssymb}
\usepackage{amsfonts}
\usepackage{amsmath}
\usepackage{epsfig}
\usepackage{color}

\allowdisplaybreaks
\makeatletter
\long\def\unmarkedfootnote#1{{\long\def\@makefntext##1{##1}\footnotetext{#1}}}
\makeatother

\newtheorem{definition}{Definition}[section]

\newtheorem{lemma}[definition]{Lemma}

\newtheorem{theorem}[definition]{Theorem}

\newtheorem{proposition}[definition]{Proposition}

\newtheorem{corollary}[definition]{Corollary}

\newtheorem{remark}[definition]{Remark}

\newtheorem{example}[definition]{Example}

\def\o{\Omega}

\def\mo{|\Omega |}
\def\m2{|\Omega | /2}
\def\M2{\frac{|\Omega |}{2}}
\def\u+{u_+^*}

\def\-p{\overline{p}}

\def\w0{{W_0^{1,p}(\Omega)}}

\def\R{\mathbb R}
\def\N{\mathbb N}

\def\rN{\mathbb R^N}
\def\rn{{{\R}^n}}
\def\T{T_{t}}

\newcommand{\hh}{{\cal H}^{n-1}}
\newcommand{\medint}{-\kern  -,395cm\int}
\newcommand{\medintinrigo}{-\kern  -,315cm\int}
\newcommand{\medelle}{-\kern  -,235cm L}
\newcommand{\medellenrigo}{-\kern  -,180cm L}
\newcommand{\qed}{\thinspace\null\nobreak\hfill
\hbox{\vbox{\kern-.2pt\hrule height.2pt
depth.2pt\kern-.2pt\kern-.2pt \hbox to1.8mm {\kern-.2pt\vrule
width.4pt \kern-.2pt\raise1.8mm\vbox to.2pt{} \lower0pt\vtop
to.2pt{}\hfil\kern-.2pt \vrule
width.4pt\kern-.2pt}\kern-.2pt\kern-.2pt \hrule height.2pt
depth.2pt \kern-.2pt}}\par\medbreak}

\title{Quasilinear elliptic problems with general growth \\ and merely integrable,  or measure, \color{black} data
} 
\numberwithin{equation}{section}

\author{
  Andrea Cianchi\\
 {\it Dipartimento di Matematica e Informatica \lq\lq U. Dini", Universit\`a di Firenze}\\ {\it Viale Morgagni 67/A, 50134 Firenze, Italy} \\{\it  e-mail: cianchi@unifi.it}
\bigskip
\\
Vladimir Maz'ya \\
  {\it   Department of Mathematics, Link\"oping University, SE-581
83 Link\"oping, Sweden}
  \\ and \\
{\it  RUDN University}\\ {\it
6 Miklukho-Maklay St, Moscow, 117198, Russia}
\\ {\it e-mail: vlmaz@mai.liu.se}
}

\date{}

\begin{document}
\maketitle
\begin{abstract} Boundary value problems for a class of quasilinear elliptic equations, with an Orlicz type growth and $L^1$ right-hand side are considered. Both Dirichlet and Neumann problems are contemplated. Existence and uniqueness of generalized solutions, as well as their regularity, are established.  The case of measure right-hand sides is also analyzed.
\color{black}
\end{abstract}

\unmarkedfootnote {
\par\noindent {\it Keywords:} Quasilinear elliptic equations,  general growth, Orlicz spaces, $L^1$-data, approximable solutions,
Dirichlet problems, Neumann problems.}
%
%\smallskip
%\par\noindent 
%This research was partly supported by the Research Project of the
%Italian Ministry of University and Research (MIUR) Prin 2012 n.2012TC7588 "Elliptic and
%parabolic partial differential equations: geometric aspects, related
%inequalities, and applications", by  GNAMPA   of the Italian INdAM (National Institute of High Mathematics), and by
% the Ministry of Education and Science of the Russian Federation, agreement n. 02.a03.21.0008.}

\section{Introduction}\label{sec1}

 We  are concerned with boundary value problems for quasilinear elliptic equations whose  right-hand side is just an integrable function,   or even a signed mesaure with finite total variation. \color{black} Problems of this kind, involving elliptic operators 
modeled upon the $p$-Laplacian,  have been systematically investigated in the literature, starting with the papers \cite{BG1, BG2}, where measure right-hand sides are also taken into account.   Contributions on this topic include \cite{ACMM, AFT, AlMe, AMST,  BeGu,  BBGGPV,  BGM, DallA, DM, R, D, DHM, Dr, DV, FS, KuM, LM, Mi0, Mi1, M1, M2, P, Pr, Ra1, Ra2}.   Uniqueness of bounded solutions  under boundary conditions prescribed outside exceptional sets had earlier been established in  \cite{Ma72} -- see also \cite[Section 15.8.4]{Mabook}.   The analysis of linear  problems goes back to  \cite{Ma1, Ma2} and   \cite{Stamp}. 
%
%Boundary value problems for quasilinear elliptic equations, modeled upon the $p$-Laplacian, with a right-hand side that is a merely integrable function (or even just a finite measure) have been systematically investigated in the literature, starting with the papers \cite{BG1, BG2}.  Linear  problems had earlier been studied in  \cite{Stamp}.  Contributions on this topic include \cite{....}.
\par %Here, we are concerned  with 
The present paper focuses on
a class of elliptic operators whose nonlinearity is not necessarily of power type. Specifically, we mainly deal  with existence, uniqueness and regularity of solutions to
Dirichlet problems of the form
%
%
%A theory of second-order derivatives for soultions to $p$-Laplacian type elliptic equations seems to be still unc
%
%
%Consider the solution $u$ to the Dirichlet problem
\begin{equation}\label{eqdirichlet}
\begin{cases}
- {\rm div} (\mathcal A(x, \nabla u))  = f(x)  & {\rm in}\,\,\, \o \\
 u =0  &
{\rm on}\,\,\,
\partial \o \,,
\end{cases}
\end{equation}
and parallel Neumann problems
\begin{equation}\label{eqneumann}
\begin{cases}
- {\rm div} (\mathcal A(x, \nabla u) ) = f(x)  & {\rm in}\,\,\, \o \\
 \displaystyle  \mathcal A(x, \nabla u) \cdot \nu =0  &
{\rm on}\,\,\,
\partial \o \,.
\end{cases}
\end{equation}
Here, and throughout the paper, $\Omega$ is a domain - i.e. a connected open set - in $\rn$,  with finite Lebesgue measure $|\Omega|$. Moreover, $\nu$ denotes the outward unit vector on $\partial \Omega$, the dot $\lq\lq \cdot "$ stands for scalar product, and 
$\mathcal A : \Omega \times \rn \to \rn$   is a Carath\'eodory function. The datum $f$ is just assumed to belong to $L^1(\Omega)$. Obviously, the additional compatibility condition
\begin{equation}\label{intf0}
\int _\o f(x)\, dx = 0\,
\end{equation}
must be imposed when the Neumann problem \eqref{eqneumann} is taken into account. Suitable regularity on $\Omega$ has also to be required in this case. 
\par  The more general problems where $f$ is replaced in \eqref{eqdirichlet} and \eqref{eqneumann}  by a signed Radon measure, with finite total variation over $\o$,  will also be discussed. \color{black}
\par
As hinted above, a critical trait of problems \eqref{eqdirichlet} and \eqref{eqneumann} is that the role played by the function $t^p$ in problems governed by power type nonlinearities
%in governing the growth in $|\xi|$ of the function $A(x, \xi)$ 
is performed by a more general Young function $B(t)$, namely a 
convex function from $[0, \infty )$ into  $[0, \infty ]$,  vanishing at $0$. Any function of this kind can obviously be written in    the form
\begin{equation}\label{B}
B(t) = \int _0^t b (\tau ) \, d\tau \qquad \hbox{for $t \geq 0$},
\end{equation}
for some  (nontrivial) non-decreasing function $b: [0, \infty) \to [0, \infty]$. 
%Functions of this kind are usually called Young functions in the literature. 
\\ Our hypotheses on the function $\mathcal A$  amount to   the  ellipiticty condition   
\begin{equation}\label{ell}
\mathcal A (x, \xi) \cdot  \xi \geq B(|\xi|) \qquad \hbox{for  a.e. $x \in \Omega$,  and for $\xi \in \rn$,}
\end{equation}
where the dot $\lq\lq \cdot "$ denotes scalar product in $\rn$, the growth condition
\begin{equation}\label{grow}
|\mathcal A (x, \xi) | \leq C( b(|\xi|) + g(x)) \qquad \hbox{for  a.e. $x \in \Omega$,  and for $\xi \in \rn$,}
\end{equation}
for some function $g \in L^{\widetilde B}(\Omega)$, where $\widetilde B$ stands for the Young conjugate of $B$, 
and the strict monotonicity condition
\begin{equation}\label{monotone}
\big(\mathcal A (x, \xi) - \mathcal A (x, \eta)\big) \cdot (\xi - \eta) >0 \quad \hbox{for  a.e. $x \in \Omega$,  and for $\xi,  \eta \in \rn$ with $\xi \neq \eta$.} 
\end{equation}
The function $B$ is assumed to be finite-valued, and to fulfill the condition
\begin{equation}\label{infsup}
1 < i_B \leq s_B < \infty,
\end{equation}
where
 \begin{equation}\label{ia}
i_B= \inf _{t>0} \frac{t b(t)}{B(t)} \qquad \hbox{and} \qquad
s_B= \sup_{t>0} \frac{t b(t)}{B(t)} \,.
\end{equation}
Assumption \eqref{infsup} is a counterpart of the assumption that $1<p<\infty$ in the case of $p$-Laplacian type equations,  in which case  $B(t)=t^p$, and $i_B=s_B=p$. In particular, the
standard $p$-Laplace operator corresponds to
the choice $\mathcal A(x, \xi)= |\xi|^{p-2}\xi$.
\par
An appropriate functional framework in the situation at hand is provided by the Orlicz-Sobolev spaces associated wtih the function $B$. Assumption \eqref{infsup} entails that these spaces are reflexive.  Hence, when the right-hand side $f$ is regular enough, a  classical theory of monotone operators applies to the Dirichlet problem \eqref{eqdirichlet} and to the Neumann problem \eqref{eqneumann}, and yields the existence and uniqueness (up to additive constants in case of the Neumann problem) of a weak solution $u$. Membership in suitable duals of the natural Orlicz-Sobolev energy  spaces associated  with problems  \eqref{eqdirichlet} and \eqref{eqneumann} is a minimal assumption on $f$ for weak solutions to be well defined, and for this theory to apply. 
\par
The point is that, as in the case of the plain $p$-Laplace equation, if $f$ is only supposed to be in $L^1(\o)$, then it need not belong to the relevant dual spaces. This is of course the interesting case. Weak solutions to problems \eqref{eqdirichlet} and \eqref{eqneumann}   are not even well defined in this case. Nevertheless, we shall show that a unique solution, in a suitably generalized sense, does exist both to the Dirichlet problem \eqref{eqdirichlet}, and, under some mild regularity assumption on $\o$, to the Neumann problem \eqref{eqneumann}. The regularity of the relevant solution will also be described. This is indeed the main purpose of our work, that also gives grounds for further regularity results under additional structure assumptions on the differential operator appearing in problems \eqref{eqdirichlet} and \eqref{eqneumann} -- see \cite{cmsecond} with this regard. 
\par The notion of solution that will be adopted is that of approximable solution, namely a function which is the  limit of a sequence of weak solutions to problems whose right-hand sides are smooth and converge to $f$.
An approximable solution  need not be even weakly differentiable. However, it is associated with a vector-valued function on $\o$, which plays the role of a substitute for its gradient in the  definition of distributional solution to \eqref{eqdirichlet} or \eqref{eqneumann}.  With  some abuse of notation, it will still be denoted by $\nabla u$.
 This notion of solution extends to the present generalized setting that introduced for operators with power type growth in \cite{BG1, DallA}.
% , and exploited in diverse contributions \cite{???}. 
\par A sharp analysis of the regularity of the approximable solution to problems \eqref{eqdirichlet} and \eqref{eqneumann} calls for the use of optimal embedding theorems for Orlicz-Sobolev spaces. In particular, the regularity of the solution to the Neumann problem \eqref{eqneumann} turns out to depend on the degree of regularity of the domain $\o$. This is   intrinsic in the problem, 
%unavoidable,  \color{black}  
as demonstrated, for instance, by the examples in \cite{Ma2}, and is related to the fact that the target spaces in the Orlicz-Sobolev embeddings in question, for functions that need not vanish on $\partial \Omega$, depend on the regularity of $\o$. 
\par  The existence and regularity results to be established for  approximable solutions to problems \eqref{eqdirichlet} and \eqref{eqneumann} continue to hold in the case when the right-hand side is a  signed measure with finite total variation over $\o$, with a completely analogous approach. On the other hand, the proof of the uniqueness of  solutions does not carry over to this more general setting. In fact, as far as we know, the problem of uniqueness of solutions to boundary value problems for  equations with a measure on the right-hand side is still open even in the basic case of differential operators with power type nonlinearities. Partial results dealing with measures that are not too much concentrated, or with differential operators fulfilling additional assumptions,  can be found in \cite{BGO, DHM, DMOP, KX}. 
\par \color{black} Let us mention that  results on  boundary value problems for quasilinear elliptic equations with Orlicz type growth and $L^1$ or measure right-hand side  are available in the literature -- see e.g. \cite{BeBe} and  subsequent papers by the same authors and their collaborators. However, different notions of generalized solutions are employed in these contributions, and, importantly,  their precise regularity properties are not discussed. Pointwise gradient bounds for local solutions to equations with Uhlenbeck  type structure and Orlicz growth, with measure-valued right-hand side, are the subject of \cite{Ba}. 
\par We  add that assumption \eqref{infsup} could probably be weakened, as to also include certain borderline situations.
%a non reflexive functional framework. 
%Methods relying upon the use of  so called complementary systems,   as, for instance, in \cite{Do, Go},  might apply in this case. 
The   results would, however, take a more technical form, and we prefer not to address this issue,  
%since it is 
not being our primary interest here.
\par  Let us conclude this section 
by recalling that elliptic systems with a right-hand side in $L^1$, or in the space of finite measures, have also been considered in the literature. The references \cite{KMsis} for local solutions, and \cite{DHM1, DHM2, DHM} for Dirichlet boundary problems,   are relevant in this connection. The results of the latter contributions are likely to admit an extension to the generalized framework of the present paper. However,  their proofs would entail substantial modifications, and the function spaces coming into play would take a somewhat different form, due to the failure of standard truncation techniques for vector-valued functions. This is of course an issue of interest, but goes beyond the scope of this work.
\color{black}

\section{Function spaces}\label{sec2}

\subsection{Spaces of measurable functions}\label{lorentzspaces}
%
%\par
%Let $\Omega$ be an open set with $|\Omega|<\infty$.
The decreasing rearrangement $u^*: [0, \infty ) \to [0, \infty ]$
of    a measurable function $u: \Omega \to \mathbb R$ is the
unique right-continuous, non-increasing function in $[0, \infty ) $
equidistributed with $u$. In formulas,  
\iffalse

 Namely, on defining the distribution
function $\mu _u: [0, \infty ) \to [0, \infty )$  of $u $ as
\begin{equation}\label{mu}
\mu _u (t) = |\{x \in \Omega :|u(x)|
>t\}| \qquad \quad \hbox{for $t\geq 0$,}
\end{equation}
%Here, $|  \cdot |$ stands for the Lebesgue measure.
 %The decreasing rearrangement $w^* : [0, \infty )
%\to [0, \infty ]$ of $u$ is  given by
the function $u^*$ is defined as

\fi
\begin{equation}\label{rearr}
u^* (s) = \inf\{ t\geq 0 : |\{|u|>t\}| \leq s\} \qquad {\rm {for}}\, \, s \in [0, \infty ).
\end{equation}
%and is the unique right-continuous non-increasing function in $[0,
%\infty ) $ equidistributed with $u$.
Clearly, $u^*(s)=0$ if $s \geq |\Omega|$. A basic property of rearrangements tells us that
\begin{equation}\label{HL}
\int_E |u(x)|\, dx \leq \int _0^{|E|} u^*(s)\, ds 
\end{equation}
 for every measurable set $E\subset \o$. 
\par
Let $\varphi : (0, |\o|) \to (0, \infty)$ be a continuous increasing function. We denote by 
$L^{\varphi, \infty}(\o)$ the Marcinkiewicz  type space of those measurable functions $u$ in $\o$ such that 
$$\sup _{s \in (0, |\o|)} \frac {u^{*}(s)}{\varphi ^{-1}(\lambda/s)}< \infty$$
for some  $\lambda >0$. Note that $L^{\varphi, \infty}(\o)$  is not always a normed space.
Special choices of the function $\varphi$ recover standard spaces of weak type. For instance, if $\varphi (t) = t ^q$ for some $q >0$, then $L^{\varphi , \infty}(\o)= L^{q, \infty}(\Omega)$, the customary weak-$L^q(\Omega)$ space. When $\varphi (t)$ behaves like $t^q (\log t)^\beta$ near infinity for some $q >0$ and $\beta \in \mathbb R$, we shall adopt the notation $L^{q, \infty}(\log L)^\beta (\o)$ for $L^{\varphi , \infty}(\o)$. The meaning of the notation $L^{q , \infty}(\log L)^\beta (\log \log L)^{-1}(\o)$ is analogous. 
%This is a rearrangement-invariant space for every $q >1$.
\par The
Orlicz spaces extend the Lebesgue spaces in the sense that the role
of powers in the definition of the norms is instead played by Young
functions, i.e. functions $B$ of the form \eqref{B}.
%A Young  function $B : [0, \infty ) \to [0, \infty ]$ is
%a left-continuous, non-trivial convex function such that $B (0)=0$. Any function $B$ of this kind admits a representation of the form
%\begin{equation}\label{intB}
%B(t) = \int _0^t b(s)\, ds \quad \hbox{for $t \geq 0$,}
%\end{equation}
%for some non-decreasing left-continuous function $b:[0, \infty ) \to [0, \infty ]$.
%If, in addition, $0< B(t) <
%\infty$ for $t>0$ and
%$$\lim _{t \to 0}\frac {B (t)}{t} =0 \qquad \hbox{and} \qquad \lim _{t \to \infty }\frac {B (t)}{t}
%=\infty,$$ then $B$ is called an $N$-function. 
The Young conjugate
of a Young function $B$ is the Young function $\widetilde B$ defined
as
\begin{equation*}%\label{young}
\widetilde B (t) = \sup \{st - B (s): s \geq 0\} \qquad \hbox{for $t
\geq 0$.}
\end{equation*}
One has that
\begin{equation}\label{intBtilde}
\widetilde B(t) = \int _0^t b^{-1}(\tau)\, d\tau \quad \hbox{for $t \geq 0$,}
\end{equation}
where $b^{-1}$ denotes the (generalized) left-continuous inverse of the function $b$ appearing in \eqref{B}.
%\\
%Note that, if  $B$ is an $N$-function, then $\widetilde B$ is an
%$N$-function as well. 
\\
As a consequence of the monotonicity of the function $b$, one has that
\begin{equation}\label{bB}
\frac t2b(t/2) \leq B(t) \leq tb(t)\qquad \hbox{for $t \geq0$.}
\end{equation}
An application of inequality \eqref{bB}, and of the same inequality  with $B$ replaced by $\widetilde B$,  yield
\begin{equation}\label{may10}
\frac 14 B(b^{-1}(t))  \leq \widetilde B (t) \leq B(2 b^{-1}(t)) \quad \hbox{for $t \geq 0$.}
\end{equation}
%
%
%Moreover,
%\begin{equation}\label{youngprop}
%s \leq B^{-1}(s) \widetilde B^{-1}(s) \leq 2s \qquad \hbox{for $s
%\geq 0$.}
%\end{equation}
A Young function (and, more generally, an increasing function) $B$
is said to belong to the class $\Delta _2$ if there exists
a constant $c$ such that
\begin{equation*}%\label{delta2}
B (2t ) \leq c B (t) \qquad \hbox{for $t\geq 0$.}
\end{equation*}
%Owing to \eqref{bB},
%\begin{equation}\label{Bbdelta2}
%\hbox{
%$B\in \Delta _2$ if and only if $b \in \Delta _2$.}
%\end{equation}
Note that, if $B\in \Delta _2$, then the first inequality in \eqref{bB} ensures that 
\begin{equation}\label{may11}
t b(t) \leq c B(t) \quad \hbox{for $t \geq 0$,}
\end{equation}
for some constant $c$, and the second inequality in \eqref{may10} implies that 
\begin{equation}\label{may12}
\widetilde B(b(t)) \leq c B(t) \quad \hbox{for $t \geq 0$,}
\end{equation}
for some constant $c>0$.
\\ The second inequality in \eqref{infsup} is equivalent to the fact that $B \in \Delta _2$; the first inequality is equivalent to $\widetilde B \in \Delta _2$. 
Under assumption  \eqref{infsup},  
\begin{equation}\label{Bmontone}
\frac{B(t)}{t^{i_B}} \quad \hbox{is inon-decreasing, \, and}\quad  \frac{B(t)}{t^{s_B}}\quad  \hbox{ is non-increasing} \quad \hbox{for $t >0$;}
\end{equation}
furthermore,
\begin{equation}\label{Btildemontone}
\frac{\widetilde B(t)}{t^{s_B'}} \quad  \hbox{ is non-decreasing, \, and} \quad \frac{\widetilde B(t)}{t^{i_B'}} \quad   \hbox{is non-increasing} \quad \hbox{for $t >0$.}
\end{equation}
Here, and in what follows, prime stands for the H\"older conjugate, namely $i_B'= \frac {i_B}{i_B-1}$ and $s_B'= \frac {s_B}{s_B-1}$. 
%Hence, given any $t_0>0$,
%\begin{equation}\label{compatti}
%\frac{B(t_0)}{t_0^{i_B}} t^{i_B} \leq B(t) \leq \frac{B(t_0)}{t_0^{s_B}}t^{s_B}\quad \hbox{for $t\geq t_0$.}
%\end{equation}
%and
%\begin{equation}\label{compattitilde}
%\frac{\widetilde B(t_0)}{t_0^{s_B'} } t^{s_B'} \leq \widetilde B(t) \leq \frac{\widetilde B(t_0)}{t_0^{i_B'} } t^{i_B'}\quad \hbox{for $t\geq t_0$.}
%\end{equation}
\\ A Young function $A$ is said to dominate another Young function $B$ near infinity  [resp. globally] if there exist  constants $c>0$ and $t_0\geq 0$ such that
\begin{equation}\label{domination}
B(t) \leq A(ct) \quad \hbox{for $t \geq t_0$ \,\, [$t \geq 0$].}
\end{equation}
The functions $A$ and $B$ are called equivalent near infinity  [globally] if they dominate each other near infinity [globally].

\par
The Orlicz space $L^{B}(\Omega )$ built upon a Young function $B$  is the Banach function space of those
real-valued  measurable functions $u$ on $\Omega$   whose Luxemburg
norm
\begin{equation*}%\label{orlicz}
\|u\|_{L^B (\Omega )} = \inf \bigg\{\lambda >0: \int _\Omega B \bigg(\frac
{|u(x)|}\lambda \bigg) \,dx \leq 1\bigg\}
\end{equation*}
is finite. The choice $B(t)=t^p$, for some $p\geq 1$, yields $L^{B}(\Omega)=L^p(\Omega)$. The Orlicz spaces associated 
%with Young functions $B(t)$ equivalent to $t^p (\log t)^\alpha$, with $p >1 $ and $\alpha \in \mathbb R$ are called Zygmund spaces, and denoted by $L^p(\log L)^\alpha (\Omega)$.
%Orlicz spaces of exponential type  are associated 
with Young functions of the form $B(t)=e^{t^\beta}-1$ for some $\beta >0$, and are denoted by $\exp L^{\beta}(\Omega)$.   The notation $\exp \exp L(\Omega)$ is adopted for Orlicz spaces associated with the Young function  $B(t)=e^{e^t}-e$.
%\\ The space $L^B(\Omega)$ is reflexive if and only if $B\in \Delta _2$ and $\widetilde B \in \Delta _2$.
%
%\begin{equation}\label{exp}
%\|f\|_{\exp L^\alpha (X, \mu)} = \inf \bigg\{\lambda >0: \frac 1{\mu (X)}\int _X e^{(\frac {|f|}\lambda)^\alpha} -1 \, d \mu \leq 1\bigg\}
%\end{equation}
%
%
\\ The H\"older type inequality
\begin{equation}\label{holder}
\int _\o |u(x) v(x)|\, dx \leq 2 \|u\|_{L^B (\o )}
\|v\|_{L^{\widetilde B} (\o )}
%
%\|v\|_{L^{\widetilde B} (\MR )} \leq {\sup}_{w \in L^B (\MR )}
%\frac{\int _\MR |v(x) w(x)|\, dx}{\|w\|_{L^B (\MR )}} \leq
%2\|v\|_{L^{\widetilde B} (\MR )}
\end{equation}
holds for every  $u \in L^B (\o )$ and $v \in L^{\widetilde B} (\o
)$.
  Given two Young functions $A$ and $B$,  one has that
\begin{equation}\label{inclusion}
L^{A}(\Omega) \to L^{B}(\Omega ) \,\,\, \hbox{if and only if} \,\,\,
\hbox{$A$ dominates $B$ near infinity.}
\end{equation}
%\par\noindent The Orlicz space $L^{B}(\MR ,
%\rN)$, with $N>1$, of $\rN$-valued measurable functions on $\MR$ is
%defined as $L^{B}(\MR , \rN)  = (L^{B}(\MR ))^N$, and is equipped
%with the norm given by $\|\v \|_{L^B(\MR , \rN)}=\|\,|\v |\,\|_{L^B(\MR)}$ for $\v \in L^B(\MR , \rN)$.
%$$\|\v \|_{L^{q, \sigma}(\MR ,
%\rN)} \approx \|\, |\v |\, \|_{L^{q, \sigma}(\MR)}\quad  \hbox{for
%every $\v \in L^{q, \sigma}(\MR , \rN)$,}$$
%   up
%  to multiplicative constants depending on  $N$.
%\begin{equation}\label{exp}
%\|f\|_{\exp L^\alpha (X, \mu)} = \inf \bigg\{\lambda >0: \frac 1{\mu (X)}\int _X e^{(\frac {|f|}\lambda)^\alpha} -1 \, d \mu \leq 1\bigg\}
%\end{equation}
 We refer to the monograph \cite{RR} for a comprehensive treatment of Young functions and Orlicz spaces.

\subsection{Spaces of weakly differentiable functions}\label{weakly}

\iffalse

Let $\Omega $ be an open set in $\rn$, with $n \geq 1$.  
Given a positive integer $m$ and a rearrangement-invariant space $X(\o)$, 
we denote by $W^mX(\o)$ the $m$-th order  Sobolev type space built upon $X(\o)$. Namely, the Banach space
\begin{multline}\label{sobolev}
W^mX(\Omega ) = \{u \in X(\Omega ):
\hbox{is $m$-times  weakly differentiable in $\Omega$} \\
\hbox{and $|\nabla ^k u| \in X(\Omega )$ for $1 \leq k
\leq m$}\}\, \end{multline}
equipped with the norm
$\|u\|_{W^mX (\Omega )} = 
\sum _{k=0}^m\|\, |\nabla ^k u\,\|_{X(\Omega )}.$
Here, $\nabla ^k u$ denotes the vector of all weak derivatives of
$u$ of order $k$. By $\nabla ^0 u$ we mean $u$. Moreover,
%Of course, by $\nabla ^0 u$ we just mean $u$, and
when $k=1$ we simply write $\nabla u$ instead of $\nabla ^1 u$.
\\ The choice $X(\o)=L^p(\o)$ reproduces the usual Sobolev space $W^{m,p}(\o)$. 
\\ When $X(\o)$ is a Lorentz space $L^{q, \sigma}(\o)$,  the  definition \ref{sobolev} yields the Lorentz-Sobolev space  $W^mL^{q, \sigma}(\o)$.  In our applications, we shall make use of the second-order Sobolev  type spaces
$W^2L^{q, \infty}(\o)$ and $W^2L^{1, \infty; \, log}(\o)$
 associated with the Marcinkievicz spaces introduced above.
\\ The Orlicz-Sobolev spaces correspond to the case when $X(\o)=L^B(\o)$ for some Young function $B$. 

\fi

Given a Young function $B$, 
we denote by $W^{1, B}(\o)$ the    Orlicz-Sobolev   space defined as 
\begin{equation}\label{sobolev}
W^{1, B}(\Omega ) = \{u \in L^B(\Omega ):
\hbox{is  weakly differentiable in $\Omega$ and $|\nabla  u| \in L^B(\Omega )$}\}\, ,
\end{equation}
equipped with the norm
$$\|u\|_{W^{1, B}(\Omega )} = \|u\|_{L^B(\Omega )} +
 \| \nabla   u \|_{L^B(\Omega )}.$$
%Here, $\nabla ^k u$ denotes the vector of all weak derivatives of
%$u$ of order $k$. By $\nabla ^0 u$ we mean $u$. Moreover,
%%Of course, by $\nabla ^0 u$ we just mean $u$, and
%when $k=1$ we simply write $\nabla u$ instead of $\nabla ^1 u$.
%
%Only first-order Orlicz-Sobolev spaces will be needed for our purposes, for which the alternate notation $W^{1,B}(\o)$ will be adopted. 
The notation $W^{1, B }_0(\Omega
)$ is employed for
%and $W^{1, B }_\bot(\Omega )$ 
the subspace of $W^{1, B
}(\Omega )$ given by
\begin{equation*}%\label{3.1.6}
W^{1, B }_0(\Omega ) = \{u \in W^{1,B}(\Omega ): \hbox{the
continuation of $u$ by $0$ outside $\Omega$ is  weakly
differentiable  in $\rn$}\}\,.
\end{equation*}
%and
%\begin{equation*}
%W^{1, B }_\bot (\Omega ) = \bigg\{u\in W^{1, B }(\Omega ):  \int
%_\Omega u(x)\,dx =0\bigg\}.
%\end{equation*}
Thanks to a Poincar\'e type inequalty in Orlicz spaces -- see \cite{Talentibound} -- the functional $\|\nabla u\|_{L^B(\o)}$ defines a norm on $W^{1, B }_0(\Omega )$ 
%and on $W^{1,B}_\bot (\o)$ 
equivalent to $\|  u\|_{W^{1,B}(\o)}$.
\\
The duals of $W^{1, B }(\Omega
)$ and  $W^{1, B }_0(\Omega
)$ will be denoted by $(W^{1, B }(\Omega
))'$ and  $(W^{1, B }_0(\Omega
))'$, respectively.
\\
 If $B \in \Delta_2$ and  $\widetilde B \in \Delta_2$, then the spaces $W^{1, B }(\Omega
)$ and  $W^{1, B }_0(\Omega
)$ 
%and $W^{1, B }_\bot(\Omega )$ 
are separable and reflexive.
\par 
%In view of the definitions of solutions to problems \eqref{eqdirichlet} and \eqref{eqneumann}, and of the proof of 
Our results involve certain  function spaces  that consist of those functions whose truncations are Orlicz-Sobolev functions. 
They extend the spaces introduced in \cite{BBGGPV} in the standard case when $B(t)=t^p$ for some $p \geq 1$.
Specifically, given any $t>0$,  let  $T_{t} : \R \rightarrow \R$ be the function
defined as
\begin{equation}\label{502}
T_{t} (s) = \begin{cases}
s & {\rm if}\,\,\, |s|\leq t \\
 t \,{\rm sign}(s) &
{\rm if}\,\,\, |s|>t \,.
%\\
% \tau &
%{\rm if}\,\,\, \tau < s  \,.
\end{cases}
\end{equation}
We set
\begin{equation}\label{503}
\mathcal T^{1,B}(\o) = \left\{u \, \hbox{is measurable in $\o$}: \hbox{$T_{t}(u)
\in W^{1,B}(\o)$ for every $t >0$} \right\}\,.
\end{equation}
The subspace   $\mathcal T^{1,B}_0(\o)$ is defined accordingly, on
replacing $W^{1,B}(\o)$  by $W_{0}^{1,B}(\o)$ on the right-hand
side of \eqref{503}. 
%The notation $\mathcal T^{1,p}(\o)$ stands for $\mathcal T^{1,B}(\o)$   when    $B(t)=t^p$ for some $p \geq 1$,  and simliarly for $\mathcal T^{1,p}_0(\o)$. 
Note that $u \in \mathcal T^{1,B}(\o)$ if and only if $u+c \in \mathcal T^{1,B}(\o)$ for every $c \in \mathbb R$.
\\ If $u \in\mathcal T^{1,B}(\o)$, then there exists
a (unique) measurable function  $Z_u : \o \to \rn$ such that
\begin{equation}\label{504}
\nabla T_{t}(u) = \chi _{\{|u|<t\}}Z_u  \qquad \quad
\hbox{a.e. in $\o$}
\end{equation}
for every $t>0$. This is a consequence of \cite[Lemma 2.1]{BBGGPV}. Here $\chi _E$ denotes
the characteristic function of the set $E$. One has that $u \in
W^{1,B}(\o)$ if and only if $u \in \mathcal T^{1,B}(\o)\cap
L^B(\o )$ and  $|Z_u| \in L^B(\o )$;  in this
case, $Z_u= \nabla u$ a.e. in $\o$. An analogous property holds provided that
$W^{1,B}(\o)$ and $\mathcal T^{1,B}(\o)$ are replaced with $W^{1,B}_0(\o)$ and $\mathcal T^{1,B}_0(\o)$, respectively . With abuse of
notation, for every $u \in \mathcal T^{1,B}(\o)$  we denote $Z_u$ simply
by $\nabla u$.
\par
A sharp embedding theorem for Orlicz-Sobolev spaces is critical in our results. As in  standard  embeddings of Sobolev type, such an embedding depends on the regularity of the  domain $\o$. The regularity of $\o$ can effectively be descirbed in terms of its isoperimetric function, or, equivalently, in terms of a relative isoperimetric inequality. 
Recall that,
given $\sigma \in [n, \infty)$,  the domain $\Omega \subset \rn$  is said to satisfy a relative isoperimetric inequality with exponent $1/{\sigma '}$ if there exists a  positive constants $C$   such that 
\begin{equation}\label{isopineq}
C \min\{|E|, |\o \setminus E|\}^{\frac 1{\sigma '}} \leq P(E;\Omega)
\end{equation}
for every measurable set $E\subset \Omega$. Here, $P(E;\Omega)$ denotes the perimeter of $E$ relative to $\Omega$, in the sense of geometric measure theory. Recall that 
$$P(E;\Omega) = \hh (\partial E \cap \Omega)$$
whenever $\partial E \cap \Omega$ is sufficiently smooth. Here, $\hh$ stands for the $(n-1)$-dimensional Hausdorff measure. The assumption that $\sigma \in [n, \infty)$ is due to the fact that inequality \eqref{isopineq} cannot hold if $\sigma < n$, whatever is $\Omega$. This can be  shown on testing inequality \eqref{isopineq} when $E$ is a ball, and letting its radius tend to $0$.
\\
We denote by     $\mathcal G _{1/{\sigma '}}$ the class of domains  in $\rn$ satisfing a relative isoperimetric inequality with exponent $1/{\sigma '}$.  These classes  were introduced in \cite{Ma0}, where membership of a domain $\o$ in  $\mathcal G _{1/{\sigma '}}$ was shown to be equivalent to    the Sobolev embedding $W^{1,1}(\o) \to L^{\sigma '}(\o)$. 
%Note that any set in $\mathcal G _{1/{\sigma '}}$ is necessarily connected. 
\\ Any bounded  Lipschitz domain belongs to $\mathcal G _{1/{n '}}$. Moreover,   the constant $C$ in \eqref{isopineq}, with $\sigma =n$, admits a lower bound depending on $\o$ only via its diameter  and its  Lipschitz constant.
\\ More generally, if $\o$ is a bounded domain in $\rn$ whose boundary is locally the graph of an H\"older continuous function with exponent $\alpha \in (0, 1]$, then
\begin{equation}\label{isopholder}
\o \in \mathcal G_{\frac {n-1}{n-1+\alpha}}.
\end{equation}
This follows on coupling   \cite[Theorem]{La} with \cite[Corollary 5.2.3]{Mabook} -- see also \cite[Theorem 1]{Crelative} for an earlier direct proof when $n=2$.
\\ Another customary class of open sets  is that of John domains. Any John domain belongs to $\mathcal G_{1/n'}$. The class of John domains is  a special instance, corresponding to the choice $\gamma =1$, of the family of the so called $\gamma$-John domains, with$\gamma \geq 1$.  The latter consists of all 
 bounded open sets $\o$ in
$\rn$ with the property that there exist a constant $c
\in (0,1)$ and a point $x_0 \in \o$ such that for every $x \in \o$
there exists a rectifiable curve $\varpi : [0, l] \to \o$,
parametrized by arclenght, such that $\varpi (0)=x$, $\varpi (l) =
x_0$, and
$${\rm dist}\, (\varpi (r) , \partial \o ) \geq c r^\gamma \qquad
\hbox{for $r \in [0, l]$.}$$
%This replacement enlarges the class in such a way that, if $\gamma >1$,  the criterion of Theorem \ref{main} can be applied  to  the Neumann problem \eqref{problemneumann}
%%\begin{equation*}
%%\begin{cases}
%% - {\rm div} (|\nabla u|^{p-2} \nabla u) =f(u) & {\rm in}\,\,\, \Omega \\
%%\frac{\partial u}{\partial {\bf n}}  =0 & {\rm on}\,\,\,
%%\partial \Omega  \,,
%%\end{cases}
%%\end{equation*}
%for certain irregular $\gamma$-John domains $\o$ for which, by contrast, Corollary \ref{mainisop} fails.
%\\ Indeed, assume that $1<p<n$ and that $\o$ is a $\gamma$-John domain with
%\begin{equation}\label{nov2016-13} 
%1\leq \gamma < \frac{p}{n-1}+1.
%\end{equation}
% Then the Poincar\'e inequality from \cite[Theorem 2.3]{KM}  and the equivalence of \eqref{sigma} and \eqref{poincare} ensure that
%\begin{equation}\label{nov2016-11}
%\nu_{\o, p} (s) \approx s^{\frac{\gamma(n-1)+1-p}{n}} \quad \hbox{for $s$ near $0$.}
%\end{equation}
%From Theorem \ref{main} we then obtain  a non trivial solution to the Neumann problem \eqref{problemneumann},  provided that
%\begin{equation}\label{65'''}
%    \lim_{t\to \pm\infty}|t|^{ 1 - \frac{np}{(\gamma (n-1)+1-p}}f(t)=0\,.
%\end{equation}
The isoperimetric inequality established in the proof of  \cite[Corollary 5]{HK}  tells us that, if $\o$ is a $\gamma$-John domain, with $1\leq \gamma < n'$,
%\begin{equation}\label{nov2016-14} 
%1\leq \gamma < n',
%\end{equation}
then
\begin{equation}\label{isopjohn}
\o \in \mathcal G_{\frac{\gamma }{n'}}.
\end{equation}
\indent Now, given $\sigma \in [n, \infty)$, let $B$ be a Young function such that
\begin{equation}\label{conv0sig}
\int _0 \bigg(\frac t{B(t)}\bigg)^{\frac 1{\sigma -1}}\,dt < \infty.
\end{equation}
Define the function
$H_\sigma : [0, \infty) \to [0, \infty)$   as 
\begin{equation}\label{Hsig}
H_\sigma (s) = \bigg(\int _0^s \bigg(\frac t{B(t)}\bigg)^{\frac 1{\sigma-1}}\,dt \bigg)^{\frac 1{\sigma '}} \qquad \hbox{for $s \geq 0$,}
\end{equation}
and let $B_\sigma $ be the Young function given by
\begin{equation}\label{Bsig}
B_\sigma (t) = B(H_\sigma ^{-1}(t)) \qquad \hbox{for $t \geq 0$,}
\end{equation}
where $H_\sigma ^{-1} : [0, \infty) \to [0, \infty]$ stands for the generalized left-continuous inverse of $H_\sigma$, and we adopt the convention that $B(\infty)=\infty$. Assume that $\Omega \in \mathcal G_{1/\sigma '}$. 
Then 
\begin{equation}\label{osnormsig}
W^{1,B}(\o) \to L^{B_\sigma}(\o)\,,
\end{equation}
where the arrow $\lq\lq \to "$ stands for continuous embedding.
Moreover, 
 $L^{B_\sigma}(\o)$ can be shown to be optimal, namely smallest possible, among all Orlicz spaces. Embedding \eqref{osnormsig} is proved in \cite{cianchisharp}  with $B_\sigma$ replaced by an equivalent Young function; the present version follows via a variant in the proof as in  \cite{cianchibound}.
Embedding \eqref{osnormsig} is equivalent to a Sobolev-Poincar\'e inequality  in 
 integral form, which asserts that
\begin{equation}\label{ossig}
\int _\o B_\sigma \Bigg(\frac{|u - {\rm med}(u)|}{c \big(\int _\o B(|\nabla  u )|)dy\big)^{1/\sigma}}\Bigg)\, dx \leq \int _\o B(|\nabla u|)dx
\end{equation}
for some constant $c=c(\Omega)$ and for every $u \in W^{1,B}(\o)$. Here, 
\begin{equation}\label{med}
 {\rm med}(u) = \inf\{t\in \mathbb R: |\{u>t\}|\leq |\Omega|/2\}\,,
\end{equation}
the median of $u$ over $\Omega$. Let us observe that the norm of embedding \eqref{osnormsig}  and the constant $c$ in inequality \eqref{ossig} admit a bound from above in terms of the constants $C$ and $\delta$ involved in inequality \eqref{isopineq}.
\\
One can verify that  $\lim _{s\to \infty
}\frac{B_\sigma^{-1}(s)}{B^{-1}(s)} =0$. Hence,  combining embedding \eqref{osnormsig} with a standard property of Orlicz spaces (e.g. as in the proof of \cite[Theorem 3]{cianchisharp}) tells us that
 \begin{equation}\label{compactsig}
 W^{1,B}(\Omega ) \to L^B(\Omega )\,.
 \end{equation}
Moreover, the Poincar\'e inequality
 \begin{equation}\label{poinc}
 \|u - {\rm med}(u)\|_{L^B(\Omega )} \leq C \|\nabla u\|_{L^B(\Omega )} 
 \end{equation}
holds for some constant $C=C(\o)$ and  for every $u \in W^{1,B}(\o)$.
% (in fact, inequality \eqref{poinc} adimts a simpler proof, via an argument analogous to that of \cite[Lemma ???]{CP}).
Notice that, If
\begin{equation}\label{divinfsig}
\int ^\infty \bigg(\frac t{B(t)}\bigg)^{\frac 1{\sigma -1}}\,dt = \infty,
\end{equation}
then the function $H_\sigma$ is classically invertible, and the fuction $B_\sigma$ is finite-valued.
On the other-hand, in the case when $B$ grows so fast near infinity that
\begin{equation}\label{convinfsig}
\int ^\infty \bigg(\frac t{B(t)}\bigg)^{\frac 1{\sigma-1}}\,dt < \infty,
\end{equation}
the function $H_\sigma^{-1}$, and hence $B_\sigma$, equals infinity for large values of its argument.
 In particular, if \eqref{convinfsig} is in force, embedding  \eqref{osnormsig} amounts to 
\begin{equation}\label{osnorminf}
W^{1,B}(\o) \to L^{\infty}(\o).
\end{equation}
Embedding \eqref{osnorminf}  can be formulated in a   form, which is suitable for our purposes. Define the function $F_\sigma : (0, \infty) \to [0, \infty)$ as 
\begin{equation}\label{Fn}
F_\sigma (t) = t^{\sigma '} \int _t^\infty \frac{\widetilde B(s)}{s^{1+\sigma'}}\, ds \quad \hbox{for $t > 0$,}
\end{equation}
and the function $G_\sigma: [0, \infty) \to [0, \infty)$ as 
\begin{equation}\label{Gn}
G_\sigma(s) = \frac s{F_\sigma^{-1}(s)} \qquad \hbox{for $s > 0$,}
\end{equation}
and $G_\sigma (0)=0$.
Then
\begin{align}\label{osinfsig}
\| u -{\rm med}(u)\|_{L^\infty(\o)} & \leq c \, G_\sigma \bigg(\int _\o B(|\nabla u|)\, dx\bigg)
\end{align}
for some constant $c=c(n)$ and for every $u \in W^{1,B}(\o)$ -- see \cite{albericocianchi}.  Let us notice that condition \eqref{convinfsig} is equivalent to the convergence of the integral in the definition of the function $F_\sigma$   \cite[Lemma 4.1]{Cianchiibero}, and that this function is increasing -- a Young function, in fact.
\\
In particular, if  $\Omega$ is a bounded Lipschitz domain, then  $\Omega \in \mathcal G_{1/n '}$, and  the results recalled above hold with $\sigma =n$. Moreover, the norm of  embedding \eqref{osnormsig}, as well as the constant $c$ in  inequalities \eqref{ossig} and  \eqref{osinfsig} admit an upper bound depending on $\Omega$ only through its diameter and its Lipschitz constant. 
\par 
The function $B_n$, given by \eqref{Bsig} with $\sigma =n$, comes into play whenever  embeddings for the space $W^{1,B}_0(\Omega)$ are in question. Actually, one has that
\begin{equation}\label{osnorm}
W^{1,B}_0(\o) \to L^{B_n}(\o).
\end{equation}
%No regularity on $\o$ is needed for \eqref{osnorm} to hold. 
Furthermore, the inequality 
\begin{equation}\label{os}
\int _\o B_n\Bigg(\frac{|u|}{c \big(\int _\o B(|\nabla  u )|)dy\big)^{1/n}}\Bigg)\, dx \leq \int _\o B(|\nabla u|)dx
\end{equation}
holds for some constant $c=c(n)$ and for every $u \in W^{1,B}_0(\o)$. Also, the embedding
 \begin{equation}\label{compact0}
 W^{1,B}_0(\Omega ) \to L^B(\Omega )
 \end{equation}
is compact.
In particular, if condition \eqref{convinfsig} holds with $\sigma =n$,
% namely 
%\begin{equation}\label{convinf}
%\int ^\infty \bigg(\frac t{B(t)}\bigg)^{\frac 1{n-1}}\,dt < \infty,
%\end{equation}
then 
\begin{align}\label{osinf}
\| u \|_{L^\infty(\o)} & \leq c \, G_n\bigg(\int _\o B(|\nabla u|)\, dx\bigg)
\end{align}
for some constant $c=c(n)$ and for every $u \in W^{1,B}_0(\o)$.
Let us stress that no regularity on $\o$ is needed for \eqref{osnorm} -- \eqref{osinf} to hold. 

\section{Main results}\label{proofexist}

Assume that $f \in L^1(\o) \cap (W^{1,B}_0(\o))'$. A function $u \in W^{1,B}_0(\o)$ is called a weak solution to the Dirichlet problem  \eqref{eqdirichlet} if 
\begin{equation}\label{231main}
\int_\o \mathcal A(x,  \nabla u) \cdot  \nabla \varphi \, dx=\int_\o f \varphi\, dx\,
\end{equation}
for every $\varphi \in W^{1,B}_0(\o)$.  In particular, the growth condition \eqref{grow} ensures that the integral on the left-hand side of \eqref{231main} is convergent. Under the present assumption on $f$, the existence  of a unique weak solution to problem \eqref{eqdirichlet} is a consequence of  the Browder-Minthy theory of monotone operators. This follows  along the same lines as in \cite[Proposition 26.12 and Corollary 26.13]{zeidler}.
%In particular, membership of the function $h$ (appearing on the right-hand side of \eqref{grow}) in $L^{\widetilde B}(\o)$ is required to guarantee that 
\par
When $f$ solely belongs to  $L^1(\o)$, the right-hand side of equation \eqref{231main} is not well defined for every $\varphi \in W^{1,B}_0(\o)$, unless condition \eqref{convinfsig} holds with $\sigma =n$, in which case, owing to \eqref{osinf}, any such function $\varphi \in L^\infty(\o)$. One might require that equation \eqref{231main} be just fulfilled  for every test function $\varphi \in C^{\infty}_0(\o)$, namely consider   distributional solutions $u$  to problem  \eqref{eqdirichlet}. Simple examples for the Laplace equation in a ball show that this kind of solution cannot be expected to belong to the natural energy spaces. Moreover, uniqueness may fail in the class of merely distributional solutions. 
Indeed, a classical example of \cite{serrin} demonstrates  that, besides  customary solutions, pathological distributional solutions may show up even in linear problems.
%
%On the other hand, requiring that equation \eqref{231} only holds for every test function $\varphi \in C^{\infty}_0(\o)$, namely that $u$ be just a distributional solution to problem  \eqref{eqdirichlet}, which need not belong to  $W^{1,B}_0(\o)$, does not guarantee its uniqueness. Indeed, a classical example of \cite{serrin} demonstrates  that, besides  pathological solutions may show up 
\par One   effective way to overcome these drawbacks  is to define a solution $u$ to  the Dirichlet problem \eqref{eqdirichlet} as a  limit of solutions to a family of approximating problems with regular right-hand sides. 

%  that are not necessarily in  $u \in W^{1,B}_0(\o)$ -- not even weakly differentiable in the worst case -- and to restore uniqueness by requiring that they be the limit of solutions to a family of approximating problems.  This leads to the following definition of approximable solution to which we alluded in Section \ref{sec1}. 
%
%
%
% and involving the function space 
%$\mathcal T ^{1,1}_0(\o)$ introduced in Section \ref{sec2}. Recall that the notation $\nabla u$ for a function $u \in \mathcal T ^{1,1}_0(\o)$ stands for the vector-valued function $Z_u$ appearing in equation \eqref{504}.

\begin{definition}\label{approxdir}
{\rm Assume that $f \in L^1(\o)$.  A measurable function $u : \o \to \mathbb R$ is called an approximable solution to the Dirichlet problem \eqref{eqdirichlet} if 
there exists a sequence $\{f_k\}\subset  L^1(\o) \cap (W^{1,B}_0(\o))'$
such that 
$f_k
\rightarrow f$ in $L^1(\o)$,
and the sequence of
weak solutions $\{u_k\}\subset W^{1, B}_0(\o)$ to problems \eqref{eqdirichlet},
with $f$ replaced by $f_k$, satisfies $u_k\rightarrow u$  a.e. in $\o$.
%\begin{equation}\label{convu}
%u_k\rightarrow u \quad
%\hbox{ a.e. in $\o$}.
%\end{equation}
}
\end{definition}

Approximable   solutions turn out to be unique. An approximable solution $u$ to \eqref{eqdirichlet}   does not belong
 to $W^{1,B}_0(\o)$ in general, and    is not necessarily weakly differentiable. However, $u$ belongs to the space  $\mathcal T ^{1,B}_0(\o)$ introduced in Section \ref{sec2}. Moreover, if $\nabla u$ is interpreted as the function $Z_u$ appearing in equation \eqref{504}, it also  fulfills the definition of 
 distributional solution to   \eqref{eqdirichlet}. These facts, and additional integrability properties of $u$ and $\nabla u$ are the content of our first main result.
The relevant integrability properties are suitably described in terms of membership in spaces of Marcinkiewicz type, defined by  the functions $\Phi_\sigma, \Psi_\sigma : (0, \infty) \to [0, \infty)$
associated with  a number $\sigma \in [n, \infty)$ and a Young function $B$ as follows.
Let 
%$$\int _0 \bigg(\frac {s}{B(s)}\bigg)^{\frac 1{n-1}}\, ds < \infty\,.$$
$\phi_\sigma  : [0, \infty) \to [0, \infty)$ be the function given by
\begin{equation}\label{phisig}
\phi _\sigma (s) = \int _0 ^s \bigg(\frac {t}{B(t)}\bigg)^{\frac 1{\sigma-1}}\, dt \quad \hbox{for $s\geq 0$,}
\end{equation}
with $B$ modified (if necessary) near zero  in such a way that \eqref{conv0sig}  holds. Then we set 
\begin{equation}\label{Phisig}
\Phi _\sigma (t) = \frac{B(\phi_\sigma ^{-1}(t))}{t} \quad \hbox{for $t>0$,}
\end{equation}
and
\begin{equation}\label{Psisig}
\Psi _\sigma (t) = \frac{B(t)}{\phi_\sigma(t)} \quad \hbox{for $t>0$.}
\end{equation}
Let us notice  that the behavior of $\Phi_\sigma$ and $\Psi_\sigma$ near infinity -- the only piece of information that will be needed on these functions -- does not depend (up to equivalence) on the behavior of $B$ near zero.

\begin{theorem}\label{existdir}
Let $\Omega$ be  any domain in $\rn$ such that  $|\o|<\infty$. Assume that the function $\mathcal A$ satisfies assumptions \eqref{ell}--\eqref{monotone}, for some finite-valued Young function $B$ fulfilling \eqref{infsup}. Let  $f \in L^1(\o)$.
Then there exists a unique approximable solution $u$ to the Dirichlet problem \eqref{eqdirichlet}. Moreover, $u \in \mathcal T^{1,B}_0(\o)$,   
\begin{equation}\label{estgraddir}
\int _\Omega b(|\nabla u|)\, dx  \leq C\int _\Omega |f|\, dx
%\int _\Omega a(|\nabla u|) |\nabla u|\, dx \leq C \int _\o |f|\, dx
\end{equation} 
for some constant  $C=C(n, |\o|, i_B, s_B)$, and
\begin{equation}\label{distrdir}
\int_\o \mathcal A(x,  \nabla u) \cdot \nabla  \varphi \, dx=\int_\o f \varphi\, dx\,
\end{equation}
for every $\varphi \in C^\infty_0(\o)$. 
   Also, 
\begin{equation}\label{umarc}
u \in \begin{cases}  L^{\Phi_n, \infty}(\o) & \quad \hbox{if 
\eqref{divinfsig} holds with  $\sigma =n$,  }
\\ 
L^\infty(\o) & \quad \hbox{otherwise,}
\end{cases}
%
%\hbox{$u \in L^{\Phi_n, \infty}(\o)$ or $u \in L^\infty(\o)$, depending on whether \eqref{divinfsig} or \eqref{convinfsig}, with $\sigma =n$, is in force,}
\end{equation}
and
\begin{equation}\label{gradmarc}
\nabla u \in L^{\Psi_n, \infty}(\o),
\end{equation}
where $\Phi_n$ and $\Psi_n$ are the functions defined as in \eqref{Phisig} and \eqref{Psisig}, with $\sigma =n$.
\\
 If $\{f_k\}$ is any sequence as in Definition \ref{approxdir}, and $\{u_k\}$ is the associated sequence of weak solutions, then
\begin{equation}\label{may1}
u_k \to u \quad \hbox{and} \quad  \nabla u_k \to \nabla u \quad \hbox{a.e. in $\Omega$},
\end{equation}
up to subsequences.
\end{theorem}

\begin{remark}\label{unique}
{\rm The uniqueness of the approximable solution to problem \eqref{eqdirichlet} ensures that, if $f \in L^1(\o) \cap (W^{1,B}_0(\o))'$, then this solution agrees with its weak solution, which is trivially also an approximable solution.}
\end{remark}

\begin{example}\label{exdir}
{\rm Assume that the function $\mathcal A$ fulfills conditions \eqref{ell}--\eqref{monotone} with $B$ satsisfying \eqref{infsup} and such  that $B(t)\approx t^p (\log t)^\beta$ near infinity, for some $p>1$ and $\beta >0$. Then, by Theorem \ref{existdir},  there exists a unique approximable solution to the Dirichlet problem \eqref{eqdirichlet}, and
\begin{equation}\label{ex1} u \in \begin{cases} L^{\frac{n(p-1)}{n-p}, \infty} (\log L)^{\frac{\beta p}{n-p}}(\o) & \qquad \hbox{if $1<p<n$,}
\\ \exp L^{\frac{n-1}{n-1-\beta}}(\o) & \qquad \hbox{if $p=n$, $\beta < n-1$,}
\\ \exp \exp L (\o) & \qquad \hbox{if $p=n$, $\beta = n-1$,}
\\ L^\infty (\Omega) & \qquad \hbox{if either $p>n$, or $p=n$ and $\beta > n-1$.}
\end{cases}
\end{equation}
Moreover, 
\begin{equation}\label{ex2}
\nabla u \in \begin{cases} L^{\frac{n(p-1)}{n-1}, \infty} (\log L)^{\frac{\beta }{n-1}}(\o) & \qquad \hbox{if $1<p<n$,}
\\ L^{n, \infty} (\log L)^{\frac{\beta n}{n-1} -1}(\o) &\qquad \hbox{if $p=n$, $\beta < n-1$,}
\\  L^{n, \infty} (\log L)^{n-1} (\log \log L)^{-1}(\o) & \qquad \hbox{if $p=n$, $\beta = n-1$,}
\\ L^{n, \infty} (\log L)^\beta (\Omega) & \qquad \hbox{if either $p>n$, or $p=n$ and $\beta> n-1$.}
\end{cases}
\end{equation}
In particular, if $\beta =0$, then the equation in \eqref{eqdirichlet} is of $p$-Laplacian type, and equation \eqref{ex1} recovers available results in the literature. Equation \eqref{ex2} also reproduces known results, with the exception  of the borderline case when $p=n$. Indeed, a theorem from \cite{DHM} yields the somewhat stronger piece of information  that $\nabla u \in L^{n,\infty}(\o)$ in this case. The proof of this result relies upon sophisticated techniques, exploiting special features of the $n$-Laplace operator, and does not seem to carry over to the present general setting.  Let us  however mention that, whereas our conclusions in \eqref{ex2} hold for any domain $\o$ with finite measure, the result of \cite{DHM} requires some regularity assumption on the $\o$. \color{black}
}
\end{example}

\begin{remark}\label{baroni}
{\rm  As mentioned is Section \ref{sec1}, pointwise gradient estimates for solutions to a class of equations  with Orlicz type growth were established in \cite{Ba} in terms of a Riesz potential type operator applied to 
an integrable, or even just measure-valued, right-hand side. The equations considered in \cite{Ba} have the special Uhlenbeck structure
\begin{equation}\label{baroni2}
\mathcal A (x, \xi) = \frac{b(|\xi|)}{|\xi|} \xi
\end{equation}
for $x \in \Omega$ and $\xi \in \rn$. Moreover, in \cite{Ba}  the first inequality in \eqref{infsup} is replaced by the stronger requirement that
\begin{equation}\label{baroni5}
2 \leq i_B\,,
\end{equation}
namely that the differential operator has a \lq\lq superquadratic" growth.
\cite[Theorem 2.1]{Ba}, combined with the boundedness of  the Riesz potential oprator of order one from $L^1(\rn)$ into $L^{n', \infty}(\rn)$, implies that 
\begin{equation}\label{baroni1}
\nabla u \in M^{\Theta , \infty}_{\rm loc}(\o)
\end{equation}
whenever the right-hand side $f \in L^1_{\rm loc}(\o)$, where 
\begin{equation}\label{baroni3}
\Theta (t) = b(t)^{n'} \quad \hbox{for $t \geq 0$.}
\end{equation}
Here, the subscript \lq\lq loc" attached to the notation of a function space on $\o$ indicates the collection of functions that belong to the relevant space on any compact subset of $\o$.
Since the function ${B(t)}/t$ is non-decreasing, one has that $\Psi _n(t) \leq \Theta (t)$ for $t \geq 0$. Thus, assertion \eqref{baroni1} is (locally) at least as strong as \eqref{gradmarc}. Properties \eqref{baroni1} and  \eqref{gradmarc} are in fact (locally) equivalent in the non-borderline case when 
\begin{equation}\label{baroni4}
s_B<n,
\end{equation}
 since, under \eqref{baroni4}, the functions $\Psi _n $ and $ \Theta$ can be shown to be bounded by each other (up to multiplicative constants), thanks to (the second assertion) in \eqref{Bmontone}. However, \eqref{baroni3} can actually yield somewhat stronger local information than  \eqref{gradmarc} if \eqref{baroni4} fails. This is the case, for instance, when $B(t)=t^n$, since \eqref{baroni3} tells us that $\nabla u \in L^{n', \infty}_{\rm loc}(\o)$, thus providing a local version of the result from \cite{DHM} recalled in Example \ref{exdir}.

\color{black}
}
\end{remark}

%
%
%\begin{theorem}\label{existneu}
%Let  $\Omega$ be any  bounded Lipschitz domain. Assume that the function $\mathcal A$ satisfies assumptions \eqref{ell}--\eqref{monotone}, for some function $B$ fulfiling \eqref{infsup}.
% Let  $f \in L^1(\o)$ be such that $\int_\o f(x)\, dx =0$. Then there exists a unique (up to additive constants) approximable solution $u$ to problem \eqref{eqneumann}, and 
%\begin{equation}\label{estgradneu}
%\int _\Omega b(|\nabla u|)\, dx  \leq C\int _\Omega |f|\, dx
%%\|b(|\nabla u|)\|_{L^1(\o)} \leq C \|f\|_{L^1(\o)} 
%%\int _\Omega a(|\nabla u|) |\nabla u|\, dx \leq C \int _\o |f|\, dx
%\end{equation} 
%for some constant $C=C(n, d_\o, L_\o, i_b, s_b)$.   Moreover, 
%\begin{equation}\label{umarcneu}
%\hbox{$u \in M^{\Phi_n}(\o)$ or $u \in L^\infty(\o)$, depending on whether \eqref{divinf} or \eqref{convinf} is in force,}
%\end{equation}
%and 
%\begin{equation}\label{gradmarcneu}
%\nabla u \in M^{\Psi_n}(\o).
%\end{equation}
%In particular, if $f_k$ is any sequence as in the definition of approximable solution, and $\{u_k\}$ is the associated sequence of (normalized) weak solutions, then
%\begin{equation}\label{may1}
%u_k \to u \quad \hbox{and} \quad  \nabla u_k \to \nabla u \quad \hbox{a.e. in $\Omega$},
%\end{equation}
%up to subsequences.
%\end{theorem}

Consider now the more general case when the function $f$ in  the equation in \eqref{eqdirichlet} is replaced by a signed Radon  measure $\mu$ with finite total variation $\|\mu\|(\o)$ on $\o$.  Approximable solutions to the corresponding Dirichlet problem 
\begin{equation}\label{eqmeas}
\begin{cases}
- {\rm div} (\mathcal A(x, \nabla u))  = \mu  & {\rm in}\,\,\, \o \\
 u =0  &
{\rm on}\,\,\,
\partial \o \,
\end{cases}
\end{equation}
can be introduced  as in Definition \ref{approxdir}, on replacing the convergence of $f_k$ to $f$ in $L^1(\o)$ by the weak convergence in measure, namely that 
\begin{equation}\label{meas1}
\lim _{k \to \infty} \int _\o \varphi f_k \, dx = \int _\o \varphi \,d\mu
\end{equation}
for every $\varphi \in C _0(\o)$. Here, $C _0(\o)$ denotes the space of continuous functions with compact support in $\o$.

\begin{theorem}\label{existdirmeas}
Let $\Omega$, $\mathcal A$ and $B$ be as in Theorem \ref{existdir}. Let  $\mu$ be a signed Radon  measure with finite total variation $\|\mu\|(\o)$.
Then there exists an approximable solution $u$ to the Dirichlet problem \eqref{eqmeas}. Moreover, $u \in \mathcal T^{1,B}_0(\o)$,   
\begin{equation}\label{estgraddir}
\int _\Omega b(|\nabla u|)\, dx  \leq C \|\mu\|(\o)
%\int _\Omega a(|\nabla u|) |\nabla u|\, dx \leq C \int _\o |f|\, dx
\end{equation} 
for some constant  $C=C(n, |\o|, i_B, s_B)$, and
\begin{equation}\label{distrdir}
\int_\o \mathcal A(x,  \nabla u) \cdot \nabla  \varphi \, dx=\int_\o  \varphi\, d\mu (x)\,
\end{equation}
for every $\varphi \in C^\infty_0(\o)$. 
   Also, the function $u$ fulfills properties \eqref{umarc} and \eqref{gradmarc}.
\end{theorem}

\color{black}

Let us next turn our attention to Neumann problems.  The definition of weak solution $u \in W^{1,B}(\o)$ to problem \eqref{eqneumann}, with $f \in L^1(\o) \cap (W^{1,B}(\o))'$ fulfilling \eqref{intf0}, reads:
\begin{equation}\label{231mainneu}
\int_\o \mathcal A(x,  \nabla u) \cdot \nabla \varphi \, dx=\int_\o f \varphi\, dx\,
\end{equation}
for every $\varphi \in W^{1,B}(\o)$. 
Analogously to the Dirichlet problem, 
the existence  and uniqueness of a  weak solution to problem \eqref{eqneumann} can be established via the theory of monotone operators.
\par
When $f\in L^1(\o)$ only, a notion of approximable solution can be introduced via a  variant of Definition \ref{approxdir}.

\begin{definition}\label{approxneu}
{\rm Assume that $f \in L^1(\o)$.  A measurable function $u$ is called an approximable solution to the Nuemann problem \eqref{eqneumann} if 
there exists a sequence $\{f_k\}\subset   L^1(\o) \cap (W^{1,B}(\o))'$ such that  $\int _\o f_k(x)\, dx =0$ for $k \in \N$,
$f_k
\rightarrow f$ in $L^1(\o)$,
and the sequence of (suitably normalized by additive constants)
weak solutions $\{u_k\}\subset W^{1, B}(\o)$ to problems \eqref{eqneumann},
with $f$ replaced by $f_k$, satisfies   $u_k\rightarrow u$  a.e. in $\o$.
%\begin{equation}\label{convu}
%u_k\rightarrow u \quad
%\hbox{ a.e. in $\o$}.
%\end{equation}
}
\end{definition}

Our existence, uniqueness, and regularity result on the Neumann problem \eqref{eqneumann} is stated in the next theorem.

\begin{theorem}\label{existneu}
Let $\Omega$ be domain in $\rn$ such that $\Omega \in \mathcal G_{1/\sigma '}$ for some $\sigma \in [n, \infty)$. Assume that the function $\mathcal A$ satisfies conditions \eqref{ell}--\eqref{monotone}, for some finite-valued Young function $B$ fulfilling \eqref{infsup}.
 Let  $f \in L^1(\o)$ be a function satisfying \eqref{intf0}. Then there exists a unique (up to additive constants) approximable solution $u$ to the Neumann problem \eqref{eqneumann}. Moreover, $u \in \mathcal T^{1,B}(\o)$,
\begin{equation}\label{estgradneu}
\int _\Omega b(|\nabla u|)\, dx  \leq C\int _\Omega |f|\, dx
%\|b(|\nabla u|)\|_{L^1(\o)} \leq C \|f\|_{L^1(\o)} 
%\int _\Omega a(|\nabla u|) |\nabla u|\, dx \leq C \int _\o |f|\, dx
\end{equation} 
for some constant $C=C(n, \Omega, i_B, s_B)$, and 
\begin{equation}\label{distrneu}
\int_\o \mathcal A(x,  \nabla u) \cdot  \nabla \varphi \, dx=\int_\o f \varphi\, dx\,
\end{equation}
for every $\varphi \in W^{1,\infty}(\o)$. 
Also, 
\begin{equation}\label{umarcneu}
u \in \begin{cases}   L^{\Phi_\sigma, \infty}(\o) & \quad \hbox{if 
\eqref{divinfsig} holds,  }
\\ 
L^\infty(\o) & \quad \hbox{otherwise,}
\end{cases}
%
%
%\hbox{$u \in L^{\Phi_\sigma, \infty}(\o)$ or $u \in L^\infty(\o)$, depending on whether \eqref{divinfsig} or \eqref{convinfsig} is in force,}
\end{equation}
and 
\begin{equation}\label{gradmarcneu}
\nabla u \in L^{\Psi_\sigma, \infty}(\o),
\end{equation}
where $\Phi_\sigma$ and $\Psi_\sigma$ are the functions defined as in \eqref{Phisig} and \eqref{Psisig}.
\\ If $\{f_k\}$ is any sequence as in Definition \ref{approxneu}, and $\{u_k\}$ is the associated sequence of (suitably normalized) weak solutions, then
\begin{equation}\label{may1}
u_k \to u \quad \hbox{and} \quad  \nabla u_k \to \nabla u \quad \hbox{a.e. in $\Omega$},
\end{equation}
up to subsequences.
\end{theorem}

A remark analogous to Remark \ref{unique} applies about the coincidence of the approximable solution with the weak solution to the Neumann problem \eqref{eqneumann} when $f \in L^1(\o) \cap (W^{1,B}(\o))'$.

\smallskip

In the case when $\o$ is a  bounded Lipschitz domain, Theorem \ref{existneu} provides us with the following piece of information.

\begin{corollary}\label{lipdomain}
Assume that $\Omega$ is a bounded  
Lipschitz domain in $\mathbb R^n$. Let $\mathcal A$, $B$ and $f$ be as in Theorem \ref{existneu}.  Then equations \eqref{umarcneu} and \eqref{gradmarcneu} hold with $\sigma=n$, and the constant $C$ in \eqref{estgradneu} admits an upper bound depending on $\Omega$ only through its
 diameter and   its Lipschitz constant.
\end{corollary}

\begin{example}\label{exneu}
{\rm  
Assume that the functions $\mathcal A$ and $B$ are as in Example \ref{exdir}. In particular,  $B(t)\approx t^p (\log t)^\beta$ near infinity. Let $\Omega \in \mathcal G_{1/\sigma '}$.  Then Theorem  \ref{existneu}  tells us that
there exists a unique (up to additive constants) approximable solution to the Neumann problem \eqref{eqneumann}, and
\begin{equation}\label{ex7} u \in \begin{cases} L^{\frac{\sigma(p-1)}{\sigma-p}, \infty} (\log L)^{\frac{\beta p}{\sigma-p}} (\o) & \qquad \hbox{if $1<p<\sigma$,}
\\ \exp L^{\frac{\sigma-1}{\sigma-1-\beta}}(\o) & \qquad \hbox{if  $p=\sigma$, $\beta < \sigma-1$,}
\\ \exp \exp L (\o)& \qquad \hbox{if  $p=\sigma$, $\beta = \sigma-1$,}
\\ L^\infty (\Omega) &  \qquad \hbox{if  either $p>\sigma$, or $p=\sigma$ and $\beta > \sigma-1$.}
\end{cases}
\end{equation}
Moreover, 
\begin{equation}\label{ex8}
\nabla u \in \begin{cases} L^{\frac{\sigma(p-1)}{\sigma-1}, \infty} (\log L)^{\frac{\beta }{\sigma-1}}(\o) & \qquad  \hbox{if  $1<p<\sigma$,}
\\ L^{\sigma, \infty} (\log L)^{\frac{\beta \sigma}{\sigma-1} -1}(\o) & \qquad \hbox{if  $p=\sigma$, $\beta < \sigma-1$,}
\\  L^{\sigma, \infty} (\log L)^{\sigma-1} (\log \log L)^{-1} (\o) &  \qquad \hbox{if  $p=\sigma$, $\beta = \sigma-1$,}
\\ L^{\sigma, \infty} (\log L)^\beta (\Omega) &  \qquad \hbox{if  either $p>\sigma$, or $p=\sigma$ and $\beta> \sigma-1$.}
\end{cases}
\end{equation}
In particular, if $\o$ is 
 a bounded  H\"older domain with exponent $\alpha \in (0, 1]$, then from  \eqref{isopholder} we infer that equations \eqref{ex7} and \eqref{ex8} hold with $\sigma = \frac{n-1+\alpha}\alpha$. 
\\ In the case when $\o$ is a $\gamma$-John domain for some $\gamma \in [1, n')$, equations  \eqref{ex7} and \eqref{ex8} follow from  \eqref{isopjohn} with $\sigma = \frac{n}{n-\gamma (n-1)}$.
}
\end{example}

We conclude this section  by mentioning that approximable solutions to the Neumann problem 
\begin{equation}\label{eqmeasneu}
\begin{cases}
- {\rm div} (\mathcal A(x, \nabla u))  = \mu  & {\rm in}\,\,\, \o \\
A(x, \nabla u) \cdot \nu  =0  &
{\rm on}\,\,\,
\partial \o \,,
\end{cases}
\end{equation}
can be  also be defined  when   $\mu$ is a signed Radon measure with  $\|\mu\|(\o)< \infty$. The counterpart of the compatibility condition \eqref{intf0} is now
\begin{equation}\label{meas0}
\mu (\o)=0\,.
\end{equation}
Definition \ref{approxneu} can be adjusted by replacing the convergence of $f_k$ to $f$ in $L^1(\o)$ by the assumption that
\begin{equation}\label{meas6}
\lim _{k \to \infty} \int _\o \varphi f_k \, dx = \int _\o \varphi \,d\mu
\end{equation}
for every $\varphi \in L^\infty (\o) \cap C(\o)$. Here, $C(\o)$ denotes the space of continuous functions in $\o$.

\begin{theorem}\label{existneumeas}
Let $\Omega$, $\mathcal A$ and $B$ be as in Theorem \ref{existneu}. Let  $\mu$ be a signed Radon  measure with finite total variation $\|\mu\|(\o)$, fulfilling \eqref{meas0}.
Then there exists an approximable solution $u$ to the Neumann problem \eqref{eqmeasneu}.  Moreover, $u \in \mathcal T^{1,B}(\o)$,
\begin{equation}\label{estgradneu}
\int _\Omega b(|\nabla u|)\, dx  \leq C \|\mu\|(\o)
%\|b(|\nabla u|)\|_{L^1(\o)} \leq C \|f\|_{L^1(\o)} 
%\int _\Omega a(|\nabla u|) |\nabla u|\, dx \leq C \int _\o |f|\, dx
\end{equation} 
for some constant $C=C(n, \Omega, i_B, s_B)$, and 
\begin{equation}\label{distrneu}
\int_\o \mathcal A(x,  \nabla u) \cdot  \nabla \varphi \, dx=\int_\o \varphi\, d\mu (x)\,
\end{equation}
for every $\varphi \in W^{1,\infty}(\o)$. 
   Also, the function $u$ fulfills properties \eqref{umarcneu} and \eqref{gradmarcneu}.
\end{theorem}

\color{black}

\section{Preliminary estimates}\label{prel}

Although we are adopting a notion of solutions different from that of \cite{BBGGPV}, our approach to Theorems  \ref{existdir} and \ref{existneu} rests upon a priori estimates for weak solutions to the problems approximating \eqref{eqdirichlet} and \eqref{eqneumann}, whose proof follows the outline of \cite{BBGGPV}. However, ad hoc Orlicz space techniques and results, such as the Sobolev type embeddings stated in Subsection \ref{weakly},  have to be exploited in the present situation. Some key steps %\color{red} tools \color{black} 
in this connection are accomplished 
%\color{red} prepared (MESSI A PUNTO; APPRONTATI) \color{black} 
in 
this section, which is devoted to  precise weak type estimates for  Orlicz-Sobolev functions satisfying a decay condition on an integral of the gradient over their level sets. This condition will be shown  to be satisfied by the approximating solutions for problems \eqref{eqdirichlet} and \eqref{eqneumann}. The ensuing estimates are crucial in the proof of the convergence of the sequence of these solutions, and of regularity properties of their limit.

\begin{lemma}\label{lemma1'} Let $\o$ be a domain in $\rn$ such that $|\Omega| < \infty$, and let $\sigma \in [n, \infty)$.  Assume that $B$ is a Young function fulfilling \eqref{conv0sig}, and let $B_\sigma$ be defined as in \eqref{Bsig}.
\par\noindent
{\rm (i}) 
%Assume that $B$ fulfils condition \eqref{conv0}. 
Let $u \in W^{1,B}_0(\o)$, and assume that there exist  constants $M>0$ and $t_0\geq 0$   such that
\begin{equation}\label{4.3}
\int_{\{|u|<t\}}B(|\nabla u|)\, dx \leq M t \qquad \hbox{for $t \geq t_0$.}
\end{equation}
\\ {\rm (a)} If   \eqref{divinfsig} holds with $\sigma =n$, then
%$$
%\int ^\infty \bigg(\frac t{B(t)}\bigg)^{\frac 1{n-1}}\,dt = \infty,
%$$
 there exists a  constant $c=c(n)$  such that
\begin{equation}\label{4.4}
|\{|u|>t\}|\leq \frac{Mt}{B_n(ct^{\frac 1{n'}}/M^{\frac 1n})} \qquad \hbox{for $t \geq t_0$.}
\end{equation}
\\ {\rm (b)} 
If    \eqref{convinfsig} holds with $\sigma =n$, then
%$$
%\int ^\infty \bigg(\frac t{B(t)}\bigg)^{\frac 1{n-1}}\,dt < \infty,
%$$
 there exists  a constant  $t_1=t_1(t_0, n, M)$ such that
\begin{equation}\label{4.4inf}
|\{|u|>t\}|=0  \qquad \hbox{for $t \geq t_1$.}
\end{equation}
\par\noindent
{\rm (ii)} Suppose, in addition, that $\o \in \mathcal G_{1/\sigma '}$.  
Let $u \in W^{1,B}(\o)$, and assume that there exist  constants $M>0$ and $t_0\geq 0$   such that
\begin{equation}\label{4.3sig}
\int_{\{|u-{\rm med}(u)|<t\}}B(|\nabla u|)\, dx \leq M t \qquad \hbox{for $t \geq t_0$.}
\end{equation}
\\ {\rm (a)} If \eqref{divinfsig} holds, 
%$$
%\int ^\infty \bigg(\frac t{B(t)}\bigg)^{\frac 1{n-1}}\,dt = \infty,
%$$
 then there exists a   constant $c=c(\o)$  such that
\begin{equation}\label{4.4sig}
|\{|u-{\rm med}(u)|>t\}|\leq \frac{Mt}{B_\sigma(ct^{\frac 1{\sigma'}}/M^{\frac 1\sigma})} \qquad \hbox{for $t \geq t_0$.}
\end{equation}
\\ {\rm (b)} 
If   \eqref{convinfsig} holds, then
%$$
%\int ^\infty \bigg(\frac t{B(t)}\bigg)^{\frac 1{n-1}}\,dt < \infty,
%$$
  there exists  a constant  $t_1=t_1(t_0, \o , M)$ such that
\begin{equation}\label{4.4infsig}
|\{|u- {\rm med}(u)|>t\}|=0  \qquad \hbox{for $t \geq t_1$.}
\end{equation}
If condition \eqref{conv0sig} is not satisfied, then an analogous statement holds, provided that $B_\sigma$ is defined as in  \eqref{Bsig}, with $B$ modified near zero in such a way that 
\eqref{conv0sig} is fulfilled. In this case, 
%the constants in inequalities \eqref{4.4}, \eqref{4.4inf}, \eqref{4.4sig} and \eqref{4.4infsig} also depend on $B$ and $|\o|$. In particular, 
the constant $c$ in \eqref{4.4} and \eqref{4.4sig} also depends on $B$. Moreover, in \eqref{4.4} and \eqref{4.4sig} the constant $t_0$ has to be replaced by another constant also depending on $B$, and $M$ has to be replaced by another constant  depending on $M$, $B$ and $|\o|$. Finally, the constant $t_1$ also depends on $B$ and $|\o|$ in \eqref{4.4inf} and \eqref{4.4infsig}.  

\end{lemma}

\medskip

\begin{remark}\label{rem3}
{\rm One can verify that that the behavior of $B_\sigma$ near infinity is independent  (up to equivalence of Young functions) of the behavior of $B$ near $0$. Thus, inequalities \eqref{4.4} and \eqref{4.4sig} are invariant (up to different choices of the involved constants, depending on the quantities described in the last part of the statement of Lemma \ref{lemma1'}) under possibly different modifications of $B$ near zero in the definition of $B_\sigma$. 
}
\end{remark}

\begin{remark}\label{remjuly}
{\rm A close inspection of the proof of Lemma \ref{lemma1'} will reveal that the same conclusions apply under the weaker assumption that $u \in \mathcal T ^{1,B}_0(\o)$ in part (i), and that $u \in \mathcal T ^{1,B}(\o)$ in part (ii). Of course, in this case $\nabla u$ has to be interpreted as the vector-valued function $Z_u$ appearing in equation  \eqref{504}.
}
\end{remark}

\begin{corollary}\label{distr} {\rm (i)} Under the assumptions of Lemma \ref{lemma1'}, part (i), for every $\varepsilon >0$, there exists $\overline t$ depending on $\varepsilon$, $M$, $t_0$, $B$, $n$ such that
\begin{equation}\label{rem1.1}
|\{|u|>t\}| < \varepsilon \quad \hbox{for $t \geq \overline t$.}
\end{equation}
{\rm (ii)} Under the assumptions of Lemma \ref{lemma1'}, part (ii), for every $\varepsilon >0$, there exists $\overline t$ depending on $\varepsilon$, $M$, $t_0$, $B$, $\o$ such that
\begin{equation}\label{rem1.1sig}
|\{|u-{\rm med}(u)|>t\}| < \varepsilon \quad \hbox{for $t \geq \overline t$.}
\end{equation}
\end{corollary}

\smallskip
\par\noindent
{\bf Proof}. Equation \eqref{rem1.1} holds trivially if the assumptions of Lemma \ref{lemma1'}, part (i), case (b) are in force. Under  the assumptions of case (a), equation \eqref{rem1.1} follows from  \eqref{4.4}. Indeed, owing to \eqref{Bmontone},  $B(t)\geq c t^{i_B}$, with $i_B>1$, for some positive constant $c$ and sufficiently large $t$ (depending on $B$), and  hence
\begin{equation}\label{limBn}
\lim _{t \to \infty} \frac{B_n(t)}{t^{n'}} = \infty.
\end{equation}
The proof of equation \eqref{rem1.1sig} is completely analogous. Now, one has to make use of the fact that the limit \eqref{limBn} also holds if $B_n$ is replaced by $B_\sigma$, and $t^{n'}$ by $t^{\sigma '}$.
\qed

\smallskip
\par\noindent
{\bf Proof of Lemma \ref{lemma1'}}. Let us focus on part (ii). Assume, for the time being, that condition \eqref{conv0sig} is fulfilled. Consider first case (a). Set
\begin{equation}\label{v}
v = u -{\rm med}(u).
\end{equation}
Clearly, $v \in W^{1,B}(\o)$, and hence $\T(v) \in W^{1,B}(\o)$ for $t>0$. Moreover, ${\rm med}(v)=0$, whence ${\rm med}(\T(v))=0$, and $\nabla v = \nabla u$ in $\Omega$.
 By  the Orlicz-Sobolev inequality \eqref{ossig} applied  to the function $\T (v)$,
\begin{equation}\label{os1}
\int _\o B_\sigma\Bigg(\frac{|\T (v)|}{C \big(\int _\o B(|\nabla \T (v)|)dy\big)^{1/\sigma}}\Bigg)\, dx \leq \int _\o B(|\nabla (\T (v))|)dx.
\end{equation}
One has that
\begin{equation}\label{os4}
\int _\o B(|\nabla \T (v)|)dx = \int _{\{|v|< t\}} B(|\nabla v|)dx \quad \hbox{for $t>0$,}
\end{equation}
and 
\begin{equation}\label{os5}
 \{|\T(v)|> t\} =  \{|v|> t\}  \quad \hbox{for $t>0$.}
\end{equation}
Thus, 
\begin{align}\label{os2}
|\{|v|>t\}| B_\sigma\bigg(\frac{ t}{C (\int _{\{|v|< t\}} B(|\nabla v|)dy)^{\frac 1\sigma} }\bigg )& \leq \int _{\{|v|>t\}} B_\sigma\Bigg(\frac{|v|}{C \big(\int _{\{|v|< t\}} B(|\nabla v|)dy\big)^{1/\sigma}}\Bigg)\, dx \\ \nonumber & \leq \int _{\{|v|< t\}} B(|\nabla v|)dx 
\end{align}
for $t>0$.  Hence, by \eqref{4.3sig},
\begin{align}\label{os3}
|\{|v|>t\}| B_\sigma\bigg(\frac{ t}{C (tM)^{\frac 1\sigma}}\bigg)  \leq M t \qquad \hbox{for $t \geq t_0$,}
\end{align}
namely \eqref{4.4sig}, with $c=C$.
\\ As far as case (b) is concerned,   an application of  inequality \eqref{osinfsig} to the function $\T (v)$, and the use of \eqref{4.3sig} yield
\begin{align}\label{os7}
\|\T (v)\|_{L^\infty(\o)} & \leq c \,G_\sigma\bigg(\int _\o B(|\nabla \T (v)|)\, dx\bigg)
\\ \nonumber & = c \,G_\sigma\bigg(\int _{\{|v|< t\}} B(|\nabla v|)\, dx\bigg) \leq c \, G_\sigma(Mt) \qquad \hbox{for $t \geq t_0$.}
\end{align}
By \eqref{os5} and \eqref{os7}, 
\begin{align}\label{os8}
|\{|v|>t\}|  =0 \qquad \hbox{if $t \geq c \,G_\sigma(Mt)$ and $t\geq t_0$.}
\end{align}
Hence, equation \eqref{4.4infsig} follows, with $t_1=  \max \big\{\tfrac{F_\sigma (cM)}{M}, t_0\big\}$, where $F_\sigma$ is deinded by \eqref{Fn}.  
\\ Assume now that condition  \eqref{conv0sig} does not hold. 
Consider the Young function $\overline B$ defined as 
\begin{equation}\label{overB}
\overline B (t) = \begin{cases}
t B(1) & \quad \hbox{if $0 \leq t \leq 1$,}
\\ B(t) & \quad \hbox{if $t>1$,}
\end{cases}
\end{equation}
and let $\overline b : [0, \infty) \to [0, \infty)$ be the left-continuous derivative of $\overline B$. Thus, 
$$\overline B (t) = \int _0^t \overline b (s)\, ds \quad \hbox{for $t\geq 0$.}$$
Clearly, $B(t) \leq \overline B(t)$ for $t\geq 0$, and condition \eqref{conv0sig} is fulfilled with $B$ replaced by $\overline B$. Let $\overline B_\sigma$ be the function defined as in \eqref{Bsig}, with $B$ replaced by $\overline B$. 
%It is easily seen that $\overline B_n$ is equivalent to $B_n$ near infinity. 
 If $u$ satisfies \eqref{4.3sig}, then 
\begin{align}\label{lem2.26}
\int_{\{|v|<t\}}\overline B(|\nabla v|)\, dx & \leq \int_{\{|v|<t, |\nabla v|>1\}}  B(|\nabla v|)\, dx + \int_{\{|v|<t, |\nabla v| \leq  1\}}\overline B(|\nabla v|)\, dx
\\ \nonumber &\leq \int_{\{|v|<t\}}  B(|\nabla v|)\, dx + B(1) |\{|v|<t\}| 
\leq t(M + B(1) |\o|) 
\end{align}
for $t \geq \max\{t_0, 1\}$.
Thereby, $u$ satisfies assumption \eqref{4.3sig}   with $B$ replaced by $\overline B$, $M$ replaced by $M+ B(1) |\o|$, and $t_0$ replaced by $\max\{t_0, 1\}$.  Hence, equation \eqref{4.4sig} holds with $B_\sigma$ replaced by $\overline B_\sigma$, $M$ replaced by $M+ B(1) |\o|$, and $t_0$ replaced by $\max\{t_0, 1\}$.
\\
By the same argument, under assumption \eqref{convinfsig}, the function $u$  still fulfills \eqref{4.4infsig}. This concludes the proof of part (ii).
\\ The proof of part (i) is completely analgous (and even simpler). One has to argue just on $u$ instead of $v$, and to make use of inequalities \eqref{os}   and \eqref{osinf} instead of    \eqref{ossig}  and \eqref{osinfsig}.
\qed

%
%
%
%
%An application of this inequality to the function $\T (u)$  tells us that
%\begin{equation}\label{os6}
%|\{|u|> t\}|=0,
%\end{equation}
%if 
%\begin{equation}\label{os7}
%t> C \|\nabla (\T (u))\|_{L^B(\o)} = C \|\nabla u\|_{L^B(\{|u|< t\})}.
%\end{equation}
%Note that here we have also made use of  equations \eqref{os5} and \eqref{os4}. On the other hand, assumption \eqref{4.3} implies that 
%\begin{equation}\label{os8}
%\int_{\{|u|<t\}}B(|\nabla u|/(Mt))\, dx \leq 1 
%\end{equation}
%if $t\geq \max\{t_0, 1/M\}$. Inequality \eqref{os8} in turn yields
%\begin{equation}\label{os9}
%\|\nabla u\|_{L^B(\{|u|< t\})}\leq Mt,
%\end{equation}
%and hence 

%\varsigma

The next lemma provides us with a uniform integrability result for functions satisfying the assumptions of Lemma \ref{lemma1'}. It will enable us to pass to the limit in the definition of distributional solution to the approximating problems.

\begin{lemma}\label{lemma2} Let $\o$ be a domain in $\rn$ such that $|\o|< \infty$. Let $B$ be a finite-valued Young function satisfying \eqref{infsup}. 
\\
{\rm (i)}
Assume that the function $u \in W^{1,B}_0(\o)$ fulfills assumption \eqref{4.3} for some  $t_0\geq 0$. Then there exists a function 
$\zeta: [0, |\o| ] \to [0, \infty )$, depending on $n$, $|\o|$, $t_0$, $B$  and on the constant $M$ appearing in \eqref{4.3}, such
that 
\begin{equation}\label{sigma0}
\lim _{s\to 0^+}\eta(s)=0,
\end{equation}
 and
\begin{equation}\label{lem2.1}
%\int _0^s b(|\nabla u|)^*(r)\, dr 
\int _E b(|\nabla u|)\, dx 
\leq \zeta (|E|) \quad \hbox{for every measurable set $E \subset \o$. }
\end{equation}
{\rm (ii)}
 Suppose, in addition, that $\o \in \mathcal G_{1/\sigma '}$ for some $\sigma \in [n, \infty)$.
Assume that the function $u \in W^{1,B}(\o)$ fulfills assumption \eqref{4.3sig} for some  $t_0\geq 0$. Then there exists a function 
$\zeta: [0, |\o| ] \to [0, \infty )$, depending on $n$, $\o$, $t_0$, $B$  and on the constant $M$ appearing in \eqref{4.3sig}, such
that \eqref{sigma0} and \eqref{lem2.1} hold.

\end{lemma}

%\begin{remark}\label{rem2}
%{\rm An inspection of the proof will reveal that the assumtpion  $u \in W^{1,B}_0(\o)$ in Lemmas \ref{lemma1} and \eqref{lemma2} can be relaxed to $u \in \mathcal T^{1,B}_0(\o)$.}
%\end{remark}

\par\noindent
{\bf Proof}. As in the proof of Lemma \ref{lemma1'}, we limit ourselves to proving part (ii). Let $v$ be defined as in \eqref{v}.
Suppose, for the time being, that 
the function $B$ satisfies condition \eqref{conv0sig}. 
\\  Assume first that condition \eqref{divinfsig} holds.  By assumption \eqref{4.3sig},
\begin{equation}\label{lem2.30}
|\{B(|\nabla v|)>s, |v|\leq t\}| \leq     \frac 1s \int_{\{B(|\nabla v|)>s, |v|\leq t\}} B(|\nabla v|)\, dx \leq \frac{Mt}s  \quad \hbox{for $t \geq t_0$ and $s > 0$.}
\end{equation}
Hence, via inequality \eqref{4.4sig}, we deduce that
\begin{align}\label{lem2.4}
|\{B(|\nabla v|)>s\}| & \leq  |\{|v|>t\}| + |\{B(|\nabla v|)>s, |v|\leq t\}| \\ \nonumber & \leq 
\frac{Mt}{B_\sigma(ct^{\frac 1{\sigma'}}/M^{\frac 1\sigma})} + \frac{M t}s \quad \hbox{for $t\geq t_0$ and $s>0$.}
\end{align}
%
%
%
%
%
%
%
%\\ Assume first that condition \eqref{divinf} holds. Let $\Lambda : [0, \infty) \times [0, \infty) \to [0, |\o|]$ be the function defined as
%\begin{equation}\label{m}
%\Lambda(t,s) = |\{x \in \o: B(|\nabla u(x)|)>s, |u(x)|>t\}| \quad \hbox{for $t \geq 0$, $s \geq 0$.}
%\end{equation}
%By inequality \eqref{os3},
%\begin{equation}\label{lem2.2}
%\Lambda(t,0) \leq \frac{Mt}{B_n(ct^{\frac 1{n'}}/M^{\frac 1n})} \qquad \hbox{for $t >t_0$.}
%\end{equation}
%Moreover, since $\Lambda$ is a decreasing function of its second variable,
%\begin{align}\label{lem2.3}
%\Lambda(0,s) & \leq \frac 1s\int _0^s \Lambda(0,r)\, dr \leq m(t,0) + \frac 1s \int _0^s (\Lambda(0,r)-\Lambda(t,r))\, dr \\ \nonumber & =\frac 1s  \int_{\{|u|<t\}}B(|\nabla u|)\, dx \leq \frac{M t}s \quad \hbox{for $t\geq t_0$ and $s>0$.}
%\end{align}
%Combining inequalities \eqref{lem2.2} and \eqref{lem2.3} yields
%\begin{equation}\label{lem2.4}
%|\{B(|\nabla u|)>s\}| \leq \frac{Mt}{B_n(ct^{\frac 1{n'}}/M^{\frac 1n})} + \frac{M t}s \quad \hbox{for $t\geq t_0$ and $s>0$.}
%\end{equation}
Set $t = \big(\tfrac 1c M^{1/\sigma} B_\sigma^{-1}(s)\big)^{\sigma '}$ in \eqref{lem2.4}, with $s \geq B_\sigma\big(ct_0^{1/\sigma '} M^{-1/\sigma}\big)$, so that $t\geq t_0$. Hence,
\begin{equation}\label{lem2.5}
|\{B(|\nabla u|)>s\}| \leq  \  \frac{2 M^{\sigma'}}{c^{\sigma'} }\frac{B_\sigma^{-1}(s)^{\sigma '}}{s} \quad \hbox{for $s \geq B_\sigma\big(ct_0 ^{1/\sigma '} M^{-1/\sigma}\big)$.}
\end{equation}
Next, 
set $\tau _0 = b\big(B^{-1}\big(B_\sigma\big(ct_0^{1/\sigma '}  M^{-1/\sigma}\big)\big)\big)$, choose $s= B(b^{-1}(\tau))$ in \eqref{lem2.5} with  $\tau \geq \tau _0$,   and make use of \eqref{Bsig} to obtain that
\begin{equation}\label{lem2.6}
|\{b(|\nabla u|)>\tau\}| \leq  \frac{2 M^{\sigma'}}{c^{\sigma'} } \frac{H_\sigma(b^{-1}(\tau ))^{\sigma '}}{B(b^{-1}(\tau ))}
%= C  M^{\sigma'}\frac{\int _0^{b^{-1}(\tau)}{\big(\frac r{B(r)}\big)^{\frac 1{\sigma-1}}}dr}{B(b^{-1}(\tau))}
%
%
%C  M^{\sigma'}\frac{\int _0^{b^{-1}(\tau)}{b(r)^{-\frac 1{\sigma-1}}}dr}{B(b^{-1}(\tau))} 
\qquad \hbox{for $\tau \geq \tau _0$.}
\end{equation}
Now, define the function $I_\sigma : (0, \infty) \to [0, \infty)$ as 
$$I_\sigma (\tau) = \frac{\tau b^{-1}(\tau )}{H_\sigma (b^{-1}(\tau))^{\sigma'}} \qquad \hbox{for $\tau >0$.}$$
Since $\frac {B(t)}t$ is  an increasing function, the function $\frac {I_\sigma(\tau)}{\tau}$ is  increasing as well. Hence, in particular, the function $I_\sigma$ is increasing.
\iffalse
\color{red} SERVE
Hence, $I_n(\tau)$ is equivalent to the Young function given by 
$$\int _0^\tau \frac {I_n(s)}{s}\,ds   \quad \hbox{for $\tau \geq 0$.}$$
\color{black}
\fi
%
%Since $s_b<\infty$, there exists a constant $c>0$ such that
%$$ct b(t) \leq B(t) \leq tb(t) \quad \hbox{for $t\geq 0$.}$$
Furthermore, by inequality \eqref{may11}, there exists a constant $c>0$  such that
\begin{equation}\label{lem2.7}
\frac{H_\sigma (b^{-1}(\tau ))^{\sigma'}}{B(b^{-1}(\tau))}I_	\sigma (\tau) \leq c \qquad \hbox{for $t \geq 0$.}
\end{equation}
Coupling \eqref{lem2.6} with \eqref{lem2.7} tells us that
\begin{equation}\label{lem2.8}
|\{b(|\nabla u|)>\tau\}| I_\sigma(\tau) \leq c M^{\sigma'} \qquad \hbox{for $\tau \geq \tau _0$,}
\end{equation}
for some constant $c>0$. Consequently, 
%\begin{equation}\label{lem2.8}
%b(|\nabla u|)^*(s)\leq I_n^{-1}\Big(\frac {C M^{n'}}s\Big) \qquad \hbox{for $s \leq \tau _0$,}
%\end{equation}
%whence
\begin{equation}\label{lem2.9}
b(|\nabla u|)^*(s)\leq I_\sigma^{-1}\Big(\frac {c M^{\sigma'}}s\Big) \qquad \hbox{if $0<s < |\{b(|\nabla u|)> \tau _0\}|$,}
\end{equation}
whence %, on setting $s_0 = \frac{CM^{n'}}{I_n(\tau _0)}$, one has that
\begin{equation}\label{lem2.10}
b(|\nabla u|)^*(s)\leq \max \bigg\{I_\sigma^{-1}\Big(\frac {c M^{\sigma'}}s\Big), \tau _0\bigg\} \qquad \hbox{for $s \in (0, |\o|)$.}
\end{equation}
From \eqref{lem2.10} and \eqref{HL}, 
inequality \eqref{lem2.1} will follow 
with $\zeta(s) = \int _0^s \max \{I_\sigma^{-1}(c M^{\sigma'}/r), \tau _0\}\, dr$,  
if we show that 
\begin{equation}\label{lem2.11}
\int _0^{|\o|} I_\sigma^{-1}\Big(\frac 1s\Big) ds 	< \infty\,.
\end{equation}
%Let $\ell \in (0, |\o|)$.  
By a change of variable and Fubini's theorem,
\begin{align}\label{lem2.12}
\int _0^{|\Omega|} I_\sigma^{-1}\Big(\frac 1s\Big) ds  = \int _{1/{|\Omega|}}^\infty \frac{I_\sigma^{-1}(t)}{t^2}\, dt 
%=  \int _{1/{|\Omega|}}^\infty \frac{1}{t^2}\int _0^{I_\sigma^{-1}(t)}dr\, dt 
 = {|\Omega|} I_\sigma^{-1}(1/{|\Omega|}) +  \int _{1/{|\Omega|}}^\infty\frac {dr}{I_\sigma(r)}.
\end{align}
Owing to Fubini's theorem again,
\begin{align}\label{lem2.13}
 \int _{1/{|\Omega|}}^\infty\frac {dr}{I_\sigma(r)} & = \int _{1/{|\Omega|}}^\infty \frac{1}{r b^{-1}(r)}\int _0^{b^{-1}(r)}\bigg(\frac s{B(s)}\bigg)^{\frac 1{\sigma-1}}ds\, dr \\ \nonumber & =\int _0^{b^{-1}(1/{|\Omega|})} \bigg(\frac s{B(s)}\bigg)^{\frac 1{\sigma-1}}ds \, \int _{1/{|\Omega|}}^\infty \frac{dr}{r b^{-1}(r)}  +
\int _{b^{-1}(1/{|\Omega|})} ^\infty \bigg(\frac s{B(s)}\bigg)^{\frac 1{\sigma-1}}\int _{b(s)}^\infty \frac{dr}{r b^{-1}(r)}\, ds.
\end{align}
Note that, by \eqref{bB} and \eqref{may11} applied with $B$ replaced by $\widetilde B$, and by \eqref{Btildemontone}, there exists a constant $c$ such that
\begin{align}\label{may14}
\int _{t}^\infty \frac{dr}{r b^{-1}(r)}  \leq \int _{t}^\infty \frac{dr}{\widetilde B(r)} 
%= \int _{t}^\infty \frac{r^{s_B'}}{\widetilde B(r)} \frac{dr}{r^{s_B'}}
 \leq 
 \frac{t^{s_B'}}{\widetilde B(t)} \int _{t}^\infty \frac{dr}{r^{s_B'}} = 
 \frac{t}{(s_B'-1)\widetilde B(t)} \leq \frac{c}{b^{-1}(t)} 
\end{align}
 for  $t>0$.
%By the second inequality in \eqref{infsup},
%\begin{align}\label{lem2.13bis}
%\int _{b(s)}^\infty \frac{dr}{r b^{-1}(r)}= \int _{s}^\infty \frac{b'(\tau )}{\tau b(\tau)}\, d\tau \leq (s_b +1) \int _{s}^\infty \frac{d\tau}{\tau ^2} = \frac{(s_b +1)}s \hbox{for $s>0$.}
%\end{align}
On the other hand, owing to \eqref{Bmontone}, one has that $B(t) \geq c t^{i_B}$ for some positive constant $c$ and sufficiently large $t$. Thus, thanks to \eqref{infsup}, the right-hand side of \eqref{lem2.13} is finite, and inequality \eqref{lem2.11} is established in this case.
\\
Assume next that condition \eqref{convinfsig} is in force. Then, by inequality \eqref{4.4infsig}, 
\begin{equation}\label{lem2.22}
|\{B(|\nabla v|)>s\}| \leq    |\{|v|>t\}| +|\{B(|\nabla v|)>s, |v|\leq t\}|  \leq    \frac{M t}s \quad \hbox{for $t\geq t_1$ and $s>0$.}
\end{equation}
%nstead of \eqref{lem2.4}. Here,   $t_1$ is the constant appearing in Lemma \ref{lemma1'}, part (ii), case (b). 
An application of inequality \eqref{lem2.22} with $t=t_1$ and $s = B(b^{-1}(\tau))$ implies that
\begin{equation}\label{lem2.23}
|\{b(|\nabla v|)>\tau\}| \widetilde B (\tau)\leq  \frac{M t_1 \widetilde B (\tau)}{B(b^{-1}(\tau))} \quad \hbox{for $\tau >0$.}
\end{equation}
%Since, by \eqref{may12}, there exists a constant $c>0$   such that
%$$\widetilde B (\tau) \leq c B(b^{-1}(\tau)) \quad \hbox{for $\tau \geq 0$,} $$
From inequalities \eqref{lem2.23} and \eqref{may12}, we deduce that
\begin{equation}\label{lem2.24}
|\{b(|\nabla v|)>\tau\}| \widetilde B (\tau)\leq  c M t_1  \quad \hbox{for $\tau >0$,}
\end{equation}
for some constant $c$.
Therefore,
\begin{equation}\label{lem2.25}
b(|\nabla v|)^* (s) \leq \widetilde B^{-1}\Big(\frac{c M t_1}s\Big)
% \max \bigg\{\widetilde B^{-1}\Big(\frac{c M t_\infty}s\Big), \tau _2\bigg\}
\quad \hbox{for $s \in (0, |\o|)$.}
\end{equation}
Inequality \eqref{lem2.1}  now follows from the convergence of the integral 
$$\int _0^{|\o|} \widetilde B^{-1}\Big(\frac{1}s\Big) \, ds,$$
which is in turn a consequence of  
%\eqref{youngprop} and 
\eqref{Btildemontone}.
%
%the fact that $\widetilde B \in \nabla _2$, and hence $\widetilde B(t) \geq C %t ^q$ for some $C>0$, $q>1$ and every $t \geq 1$.
\\ It remains to remove the temporary assumption \eqref{conv0sig}. If \eqref{conv0sig} fails, one can argue as above, with $B$ and $B_\sigma$ replaced with the functions $\overline B$ and $\overline B _\sigma$ defined as in the proof of Lemma \ref{lemma1'}, and the constants $M$, $c$, $t_0$ and $t_1$ appearing in inequalities \eqref{4.3sig}--\eqref{4.4infsig} modified accordingly.  \qed

\section{Proof of the main results}\label{mainproof}

We begin with an estimate of the form \eqref{estgraddir} for the weak solution to the Dirichlet problem \eqref{eqdirichlet} with regular right-hand side $f$, and a parallel estimate for the weak solution to the Neumann problem \eqref{eqneumann}. They follow from a slight extensions of  results of \cite{talenti} and \cite{cimanifolds}, respectively (see also of \cite{Ma2} for linear equations).

\begin{proposition}\label{talenti} Assume that conditions  \eqref{ell}--\eqref{infsup}  are fulfilled.  Let $\o$ be a domain in $\rn$ with $|\o|<\infty$.
\\ 
{\rm (i)}  Assume that $f \in L^1(\o)\cap (W^{1,B}_0(\o))'$. Let $u$ be the weak solution to the Dirichlet problem \eqref{eqdirichlet}. Then there exists a constant $C=C(n, i_B, s_B)$ such that
\begin{equation}\label{talentigrad}
\int _\Omega b(|\nabla u|)\, dx  \leq C |\o|^{\frac 1n} \int _\Omega |f|\, dx.
\end{equation}
\\ 
{\rm (ii)}  Suppose, in addition, that $\Omega \in \mathcal G_{1/\sigma '}$ for some $\sigma \in [n, \infty)$. Assume that $f \in L^1(\o)\cap (W^{1,B}(\o))'$ and satifies condition \eqref{intf0} . Let $u$ be a weak solution to the Neumann problem \eqref{eqneumann}. Then there exists a constant $C=C(n, \o, i_B, s_B)$ such that
\begin{equation}\label{neumanngrad}
\int _\Omega b(|\nabla u|)\, dx  \leq C  \int _\Omega |f|\, dx.
\end{equation}
In particular, if $\Omega$ is a bounded  Lipschitz domain, then the constant $C$ in \eqref{neumanngrad} depends only $\o$ only through its diameter and its Lipschitz constant of $\o$.
\end{proposition}
{\bf Proof, sketched}.  Part (i). If the function $B(b^{-1}(t))$ is convex, then an application of  \cite[Theorem 1, part (vi)]{talenti}, with the choice $M(t)=b(t)$ (in the notation of that theorem) entails that there exists a constant $C=C(n)$ such that
\begin{equation}\label{tal1}
\int _\o b(|\nabla u|)\, dx \leq C \int _0^{|\o|} \frac 1{s^{1/n'}}\int _0^s f^*(r)\,dr \, ds.
% \leq C \|f\|_{L^1(\o)} \int _0^{|\o|}\frac 1{s^{1/n'}}\, ds = nC |\o|^{1/n} \|f\|_{L^1(\o)} .
\end{equation}
Hence, inequality \eqref{talentigrad} follows, since
$$\int _0^s f^*(r)\,dr \leq  \|f\|_{L^1(\o)}  \qquad \hbox{for $s \in [0, |\o|]$.}$$
In general, the function $B(b^{-1}(t))$ is only equivalent to a convex function. Indeed, 
owing to inequalities \eqref{may10} and \eqref{may12}, there exist constants $c_1$ and $c_2$, depending on $i_B$ and $s_B$, such that
\begin{equation}\label{Bb-1}
c_1\widetilde B(t) \leq B(b^{-1}(t)) \leq c_2 \widetilde B(t) \quad \hbox{for $t>0$.}
\end{equation}
Following the lines of the proof of \cite[Theorem 1, part (vi)]{talenti}, with the choice   $K(t) = \widetilde B(t)$ (again in the notation of that proof), and making use of 
 \eqref{Bb-1},  tell us that inequality \eqref{tal1} continues to hold, with $C$ depending also on $i_B$ and $s_B$. Thus, inequality  \eqref{talentigrad}  holds also in this case.
\\
Part (ii) The proof of  inequality \eqref{neumanngrad} relies upon a counterpart of  the result of \cite{talenti} for the Neumann problem \eqref{eqneumann}, contained in \cite[Theorem 3.1]{cimanifolds}.  A variant of the proof of that theorem as hinted above, exploiting \eqref{Bb-1} again, implies that
\begin{equation}\label{tal2}
\int _\o b(|\nabla u|)\, dx \leq C \int _0^{|\o|} \frac 1{s^{1/\sigma '}}\int _0^s f^*(r)\,dr \, ds \leq C \|f\|_{L^1(\o)} \int _0^{|\o|}\frac 1{s^{1/\sigma'}}\, ds = \sigma C |\o|^{1/\sigma} \|f\|_{L^1(\o)}\,,
\end{equation}
where the constant $C$ also depends on the constant appearing in the relative isoperimetric inequality \eqref{isopineq}. This yields \eqref{neumanngrad}. 
\\ The assertion about bounded  Lipschitz domains rests upon the fact that, for these domains, inequality \eqref{isopineq} holds,  with $\sigma =n$,  and a constant $C$ depending on $\o$ only via its diameter  and its Lipschitz constant.
\qed
%
%there exists a constant $C=C(n, i_b, s_b)$ such that
%$$
%\int _\o b(|\nabla u|)\, dx \leq C \int _0^{|\o|} \frac 1{s^{1/n'}}\int _0^s f^*(r)\,dr \, ds \leq C \|f\|_{L^1(\o)} \int _0^{|\o|}\frac 1{s^{1/n'}}\, ds = nC |\o|^{1/n} \|f\|_{L^1(\o)} .$$
%Hence, inequality  \eqref{talentigrad}  follows. \qed

%
%
%\section{Proof of Theorems \ref{existdir} and \ref{existneu}}\label{proofsmain}
%

We are now ready to prove  our main results. We limit ourselves to provide a full proof of Theorems \ref{existneu}, and a sketch of the (completely analogous) proof of Theorem \ref{existneumeas}. The proofs of Theorems \ref{existdir} and \ref{existdirmeas} follow along the same lines, with some slight simplification.  One has basically to replace the use of parts (ii) of the preliminary results established above and in Section \ref{prel} with the corresponding parts (i), whenever they come into play. The details will be omitted, for brevity.

\medskip
\par\noindent
{\bf Proof of Theorem \ref{existneu}}.  For ease of presentation, we split this proof in steps. 
\\ {\bf Step 1}  \emph{There exists a sequence of problems,  approximating \eqref{eqneumann},  with regular right-hand sides, and a corresponding sequence $\{u_k\}$ of weak solutions.}
\\
Let  $\{ f_k \}\subset C^\infty _0(\o)$ be a sequence  of functions
such that 
$\int _\o f_k\, dx = 0$
for $k \in \mathbb N$, and
\begin{equation}\label{404}
f_k \rightarrow f \qquad \hbox{in }L^1(\o).
 % \qquad \hbox{as }
%k\rightarrow \infty
\end{equation}
%and
%\begin{equation}\label{405}
%\int_\o f_k \,dx  =0 \qquad \hbox{for } k\in \N.
%\end{equation}
%Indeed, if $1\le q<\infty$, any sequence $\{ f_k \}$ of continuous
%compactly supported functions fulfilling \eqref{404} and
%\eqref{405} does the job; when $q=\infty$, it suffices to take
%$f_k=f$ for $k\in \N$, since $V^{1,p}(\o) \to L^1(\o)$ provided
%that \eqref{402} is in force, by \eqref{2002'}. 
We may also
clearly assume that
\begin{equation}\label{406}
\|f_k\|_{L^1(\o)}\le 2\|f\|_{L^1(\o)}, \qquad \hbox{for } k\in \N.
\end{equation}
%ssumption \eqref{infsup} is equivalent to the fact that $B, \widetilde B \in \Delta _2$. Hence, $W^{1,B}_0(\Omega)$ is reflexive. As a consequence of standard results on monotone operators (see e.g. \cite{browder}), 
%
%By \cite[???]{???},
%%
%%\cite[Theorem 2.13]{cm_bound}, 
As observed in Section \ref{proofexist}, under our assumptions on $f_k$
for each $k\in \N$ there exists a
unique weak solution $u_k\in W^{1,B}(\o)$ to the Neumann problem
\begin{equation}\label{407}
\begin{cases}
- {\rm div} (\mathcal A(x, \nabla u_k)) = f_k(x)  & {\rm in}\,\,\, \o \\
\displaystyle  \mathcal A(x, \nabla u_k) \cdot  \nu =0 &
{\rm on}\,\,\,
\partial \o \,,
\end{cases}
\end{equation}
normalized in such a way that
\begin{equation}\label{408}
{\rm med } (u_k)=0.
\end{equation}
Hence,
\begin{equation}\label{408bis}
\int_\o \mathcal A(x, \nabla u_k) \cdot \nabla \varphi \, dx=\int_\o f_k \, \varphi\,
dx\, 
\end{equation}
for every $\varphi \in W^{1, B}(\o)$.
%
%We split the proof of the existence of an approximable solution to
%\eqref{1} in steps. The outline of the argument is related to that
%of \cite{BBGGPV, DMOP}. \par\noindent

\medskip
\noindent{\bf Step 2.} \emph{There exists a measurable function
$u: \Omega \to \R$ such that}
\begin{equation}\label{408ter}
u_k\rightarrow u \qquad \hbox{a.e. in }\o,
\end{equation}
\emph{up to subsequences. Hence, $u$ is an approximable solution to problem \eqref{eqneumann}.
% property (ii) of the definition
%of approximable solution holds
}
\\
Given any $t, \tau>0$, one has that
\begin{equation}\label{409}
|\{|u_k-u_m|>\tau  \}|\le |\{| u_k|>t \}|+|\{| u_m|>t \}|+
|\{|\T(u_k)-\T(u_m)|>\tau  \}|,
\end{equation}
for $k,m\in \N$. 
On choosing $\varphi=\T(u_k)$ in \eqref{408bis} and making use of assumptions \eqref{ell} and
\eqref{406} one obtains that
\begin{align}\label{414}
\int_\o B(|\nabla \T(u_k)|)\, dx & =\int_{\{|u_k|<t\}}B(|\nabla u_k|)\, dx  \le \int_{\{|u_k|<t\}}\mathcal A(x, \nabla u_k)  \cdot
\nabla u_k\, dx
\\ \nonumber  & = 
 \int_{\o}\mathcal A(x, \nabla u_k)  \cdot
\nabla \T(u_k)\, dx
=
\int_{\o} f_k \, \T(u_k)\, dx
\le 2t\|f\|_{L^1(\o)},
\end{align}
for $k\in \N$.
By \eqref{414} and Corollary \ref{distr}, part (ii),
%
%
%
%By \eqref{406} and \cite[???]{??}, the sequence $\{u_k\}$ is uniformly bounded in $L^1(\o)$. 
%
%
%
%\eqref{A} and \eqref{406},
%\begin{equation}\label{410}
%(u_k)_\pm ^* (s) \leq
%\nu _p(s)^{\frac{1}{1-p}}\|(f_k)_\pm\|_{L^1(\o)}^{\frac{1}{p-1}}\le
%2^{\frac{1}{p-1}}|\o|^{\frac{1}{q'(p-1)}}\nu _p(s)^{\frac{1}{1-p}}\|f\|_{L^q(\o)}^{\frac{1}{p-1}}
%\quad \hbox{for }s\in(0, |\o |/2),
%\end{equation}
%and for $k\in \N$, whence
%\begin{equation}\label{411}
%\mu_{(u_k)_\pm}(t)\le
%\nu _p^{-1}\left(\frac{2|\o|^{\frac{1}{q'}}\|f\|_{L^q(\o)}}{t^{p-1}}\right),
%\qquad \hbox{for }t>0,
%\end{equation}
%and for $k\in \N$. Here, $\nu _p^{-1}$ denotes the generalized
%left-continuous inverse of $\nu _p$. 
fixed any $\varepsilon>0$,
the number $t$ can be chosen so large that
\begin{equation}\label{412}
 |\{| u_k|>t \}|<\varepsilon  \quad \hbox{and} \quad  |\{| u_m|>t \}|<\varepsilon
\end{equation}
%Next, fix any smooth open set $\o_\varepsilon\subset \subset \o$ such that
%\begin{equation}\label{413}
%|\o \setminus \o_\varepsilon|<\varepsilon.
%\end{equation}
%On choosing $\varphi=\T(u_k)$ in \eqref{408bis} and making use of
%\eqref{406} we obtain that
%\begin{equation}\label{414}
%\int_\oB(|\nabla \T(u_k)|)\, dx =\int_{\{|u_k|<t\}}B(|\nabla u_k|)\, dx \le \int_{\{|u_k|<t\}}a(\nabla u_k)\nabla u_k \cdot
%\nabla u_k\, dx\le 2t\|f\|_{L^q(\o)},
%\end{equation}
%for $k\in \N$. 
for every $k, m \in \N$. Thanks to condition \eqref{408}, one has that ${\rm med}(\T (u_k))=0$. Thus, owing to inequalities \eqref{poinc}
and \eqref{414},   
the sequence $\{\T(u_k)\}$  is
bounded in $W^{1,B}(\o)$. By the compactness of embedding \eqref{compact0},   the sequence $\{\T(u_k)\}$
converges  (up to subsequences) to some function in
$L^B(\o)$. In particular, $\{\T(u_k)\}$ is a Cauchy
sequence in measure in $\o$. Thus,
\begin{equation}\label{416}
|\{ |\T(u_k)-\T(u_m)|>\tau  \}|\le  \varepsilon
\end{equation}
provided that $k$ and $m$ are sufficiently large. By \eqref{409},
\eqref{412} and \eqref{416}, $\{ u_k\}$ is (up to subsequences) a
Cauchy sequence in measure in $\o$, and hence there exists a
measurable function $u: \o\rightarrow \R$ such that \eqref{408ter}
holds.
\medskip

\noindent{\bf Step 3.}
\begin{equation}\label{416bis}
\{\nabla u_k\} \hbox{\emph{ is a Cauchy sequence in measure}.}
\end{equation}
Let $t>0$. Fix any $\varepsilon >0$. Given any $\tau,\delta>0$, we have that
\begin{align}\label{417}
|\{|\nabla u_k-\nabla u_m|>t \}|& \le |\{|\nabla  u_k|>\tau
\}|+|\{| \nabla u_m|>\tau\}|+ |\{|u_k-u_m|>\delta  \}| \\
&+|\{|u_k-u_m|\le \delta,  |\nabla  u_k|\le \tau , |\nabla
u_m|\le \tau , |\nabla u_k-\nabla u_m|>t \}|,\nonumber
\end{align}
for $k,m\in \N$. By \eqref{neumanngrad} and \eqref{406},
\begin{equation}\label{417bis}
\int _\o b(|\nabla u_k|)\, dx \leq  C   \|f\|_{L^1(\o
)}\,,
\end{equation}
for some constant $C=C(n, \o, i_B, s_B)$. Hence,
\begin{equation}\label{418}
 |\{|\nabla  u_k|>\tau \}|\le b^{-1}\Big( \frac C \tau \|f\|_{L^1(\o )} \Big),
\end{equation}
for $k\in \N$ and for some constant $C$ independent of $k$. Thus
$\tau$ can be chosen so large that
\begin{equation}\label{419}
 |\{|\nabla  u_k|>\tau \}|< \varepsilon \qquad \hbox{for }k\in\N.
\end{equation}
For such a choice of $\tau$,  set
\begin{equation}\label{420}
G=\{ |u_k-u_m|\le \delta, |\nabla u_k|\le \tau,  |\nabla u_m|\le
\tau, |\nabla u_k -\nabla u_m|\ge t\}.
\end{equation}
We claim that
%, if \eqref{419} is fulfilled, then
there exists
$\delta >0$  such that
\begin{equation}\label{421}
|G|< \varepsilon .
\end{equation}
To verify our claim, observe that, if we define
$$
S=\{(\xi, \eta)\in {\R} ^{n} \times \rn: |\xi|\le \tau, |\eta|\le \tau,
|\xi-\eta|\ge t \},
$$
and $\psi :\o \to [0,\infty)$ as
$$
\psi (x)=\inf\big\{\big(\mathcal A(x, \xi)-\mathcal A (x, \eta) \big)\cdot (\xi-\eta): (\xi,
\eta)\in S\big\},
$$
then $\psi (x)\ge 0$ and
\begin{equation}\label{422}
|\{\psi (x)=0\}|=0.
\end{equation}
This is a consequence of inequality \eqref{monotone} and of the fact that the set
$S$ is compact and the function $\mathcal A(x, \xi)$ is continuous in $\xi$ for every
$x$ outside a subset of $\o$ of Lebesgue measure zero.
%
%
%
%$$
%l=\inf\{[a(|\xi|)\xi -a(|\eta|)\eta ]\cdot (\xi-\eta):  (\xi,
%\eta)\in S\},
%$$
%then  
%\begin{equation}\label{422}
%l >0.
%\end{equation}
%This is a consequence of inequality \eqref{monotone}, and of the fact that
%$S$ is compact and $a(|\xi|)$ is continuous in $\xi$.
Therefore,
\begin{align}\label{423}
\int _G \psi (x) \, dx & \le\int_G \big(\mathcal A(x, \nabla u_k) -\mathcal A(x, \nabla u_m)\big)\cdot (\nabla u_k-\nabla u_m)\, dx\\ \nonumber
&\le \int_{\{|u_k-u_m|\le \delta\}}\big(\mathcal A(x, \nabla u_k) -\mathcal A(x, \nabla u_m)\big)\cdot (\nabla u_k-\nabla u_m)\, dx\\ \nonumber
&= \int_\o\big(\mathcal A(x, \nabla u_k) -\mathcal A(x, \nabla u_m)\big)\cdot \nabla(T_{\delta}(u_k - u_m)) \, dx \\ \nonumber
&= \int_\o(f_k-f_m)T_{\delta}(u_k - u_m)dx\le 4
\delta \|f\|_{L^1(\o)},
\end{align}
where the last equality follows on making use of $T_{\delta}(u_k -
u_m)$ as test function in \eqref{408bis},    and in the same equation with $k$ replaced with $m$, and 
subtracting the resulting equations. 
Inequality \eqref{421} is a consequence of \eqref{423}. \par\noindent Finally,  since, by
Step 1, $\{u_k\}$ is a Cauchy sequence in measure in $\o$,
\begin{equation}\label{424}
|\{ |u_k-u_m|>\delta\}|<\varepsilon\,,
\end{equation}
if $k$ and $m$ are sufficiently large. Combining \eqref{417}, \eqref{419}, \eqref{421} and \eqref{424} yields
$$
|\{ |\nabla u_k-\nabla u_m|>t\}|<4\varepsilon\,,
$$
for sufficiently large $k$ and $m$. Owing to the arbitrariness of $t$, property \eqref{416bis} is thus established.
\medskip

\noindent{\bf Step 4.} $u\in \mathcal T^{1,B}(\o)$, \emph{and}
\begin{equation}\label{425}
\nabla u_k \rightarrow \nabla u \qquad \hbox{a.e. in }\o\,,
\end{equation}
\emph{up to subsequences. Here, $\nabla u$ is the generalized
gradient of $u$ in the sense of the function $Z_u$ appearing in \eqref{504}. }
%$W^{1,p}_T(\o)$.
\\
Since $\{\nabla u_k\}$ is a Cauchy sequence in measure, there exists a measurable function $W : \o\rightarrow \rn$ such that
\begin{equation}\label{426}
\nabla u_k \rightarrow W \qquad \hbox{a.e. in }\o
\end{equation}
(up to subsequences). As observed in Step 2,
%
%Fix any $t>0$. Owing to \eqref{408}, ${\rm med}(\T(u_k))=0$ for $k \in \mathbb N$. Hence, by \eqref{414} and \eqref{poinc}, 
the sequence $\{\T(u_k)\}$
is bounded in $W^{1,B}(\o)$. Since   the Sobolev space $W^{1,B}(\o)$ is reflexive,  there exists a function
$\overline u_t\in W^{1,B}(\o)$ such that
\begin{equation}\label{427}
\T(u_k)  \rightharpoonup  \overline u_t \qquad \hbox{weakly in
}W^{1,B}(\o)
\end{equation}
(up to subsequences). By Step 2, $\T(u_k)\rightarrow \T(u)$ a.e.
in $\o$, and hence
\begin{equation}\label{428}
  \overline u_t =\T (u) \qquad \hbox{a.e. in }\o.
  \end{equation}
Thereby, $\T (u)  \in W^{1,B}(\o)$, and
\begin{equation}\label{429}
\T(u_k)  \rightharpoonup \T (u) \qquad \hbox{weakly in
}W^{1,B}(\o).
\end{equation}
By the arbitrariness of $t$, one has that 
 $u\in \mathcal T^{1,B}(\o)$, and
\begin{equation}\label{430}
\nabla \T(u) =\chi_{\{|u|<t\}}\nabla u  \qquad \hbox{a.e. in
}\o\,
\end{equation}
for $t>0$.
%where $\nabla u$ is the generalized gradient in the sense of $W^{1,p}_T(\o)$.
Owing to \eqref{408ter} and \eqref{426},
$$
\nabla \T(u_k) =\chi_{\{|u_k|<t\}} \nabla u_k \rightarrow
\chi_{\{|u|<t\}}W \qquad \hbox{a.e. in }\o\,
$$
for $t>0$.
Hence, by \eqref{429}
\begin{equation}\label{431}
\nabla \T(u) =\chi_{\{|u|<t\}}W  \qquad \hbox{a.e. in }\o\,
\end{equation}
for $t>0$.
Coupling \eqref{430} with
\eqref{431} yields
\begin{equation}\label{432}
 W =\nabla u\qquad \hbox{a.e. in }\o.
\end{equation}
Equation \eqref{425} is a consequence of \eqref{426} and \eqref{432}.
\smallskip
\par\noindent
{\bf Step 5.} \emph{$u$ fulfills inequality \eqref{estgradneu} and equation \eqref{distrneu}}.
\par\noindent
From \eqref{425} and \eqref{417bis}, via Fatou's lemma, we deduce
that
$$
\int _\o b(|\nabla u|)\, dx \le C\|f\|_{L^1(\o)},
$$
%for some constant $C$ independent of $f$.
%depending on
%$\sigma, p$ and $|\o|$.
namely \eqref{estgradneu}.
\par\noindent
As far as equation \eqref{distrneu} is concerned, observe that, by \eqref{425},
\begin{equation}\label{433}
\mathcal A(x, \nabla u_k )\rightarrow \mathcal A(x,\nabla u) \qquad \hbox{for a.e.
$x \in \o$ .}
\end{equation}
Now, fix any function $\varphi \in W^{1,\infty}(\o)$ and any measurable set
$E\subset \o$. Owing to inequality \eqref{grow} and Lemma \ref{lemma2}, part (ii) (whose assumptions are fulfilled  thanks to \eqref{414}),
\begin{align}\label{434}
\int_E |\mathcal A(x, \nabla u_k) \cdot \nabla \varphi|\, dx&\le  C \|\nabla
\varphi\|_{L^\infty(\o)}\bigg( \int_E b\big(| \nabla u_k|\big) \,
dx + \int _E g(x) \, dx\bigg)
\\ \nonumber &\le C\|\nabla \varphi  \|_{L^\infty(\o)} \bigg(\zeta
(|E|) + \int _E g(x) \, dx\bigg)
%
%\nonumbe
% \\
%&\le C \|\nabla
%\phi\|_{L^\infty(\o)}\left(\|f\|_{L^q(\o)}|F|^{\frac{\sigma-p+1}{\sigma}}+\|h\|_{L^{p'}(\o)}|F|^{\frac{1}{p}}
% \right)\nonumber
\end{align}
for some function $\zeta: [0, |\o|  \to [0, \infty )$ such
that $\lim _{s\to 0^+}\zeta(s)=0$. From \eqref{433} and
\eqref{434}, via Vitali's convergence theorem, we deduce that the
left-hand side of \eqref{408bis} converges to the left-hand side
of \eqref{distrneu} as $k\rightarrow\infty$. The right-hand side of
\eqref{408bis} trivially converges to the right-hand side of
\eqref{distrneu}, by \eqref{404}. Hence, equation \eqref{distrneu} follows. 
% This completes the proof of the
%present step, and hence also the proof
%of the existence of an approximable solution to \eqref{eqdirichlet}.

\smallskip
\par\noindent{\bf Step 6.} \emph{The solution $u$ is unique (up to additive constants).}
\par\noindent
 Assume that $u$ and $\overline u$ are approximable
solutions to problem \eqref{eqdirichlet}.  Then, there exist sequences
$\{f_k\}$ and $\{\overline f_k\}$ in $L^1(\o) \cap (W^{1,B}(\o))'$ 
%\subset L^{q}(\o)\cap (V^{1,p}(\o
%))'$ 
such that  $\int _\o f_k\, dx = \int _\o \overline f _k\, dx =0$ for $k \in \N$, 
$f_k\rightarrow f$ and $\overline
f_k\rightarrow f$ in $L^1(\o)$,   the weak solutions $u_k$ to
problem \eqref{407} fulfill $u_k \rightarrow u$ a.e. in $\o$,  and
 the weak solutions
 $\overline u_k$ to problem \eqref{407} with $f_k$ replaced by $\overline f_k$
fulfill  $\overline u_k \rightarrow
\overline u$ a.e. in $\o$. Fix any $t>0$, and choose the test
function $\varphi = \T(u_k-\overline u_k)$ in \eqref{408bis}, and in
the same equation with $u_k$ and $f_k$ replaced by $\overline u_k$
and $\overline f_k$, respectively. Subtracting the resulting
equations yields
\begin{equation}\label{435}
\int_\o\chi_{\{|u_k-\overline u_k|\le t\}} \big(\mathcal A(x, \nabla u_k) -\mathcal A(x, \nabla \overline u_k )\big)\cdot (\nabla u_k-\nabla \overline u_k)\,
dx=\int_\o (f_k-\overline f_k)\T(u_k-\overline u_k)\, dx
\end{equation}
for $k\in \N$. Since $|\T(u_k-\overline u_k)|\le t$ in $\o$ and
$f_k-\overline f_k\rightarrow 0$ in $L^1(\o)$, the right-hand side
of \eqref{435} converges to $0$ as $k\rightarrow\infty$. On the
other hand, the arguments of Steps 3 and 4
%analogous to those exploited above in the
%proof of the existence 
tell us that $\nabla u_k \rightarrow \nabla
u$ and $\nabla \overline u_k \rightarrow \nabla \overline u$ a.e.
in $\o$ (up to subsequences).  Hence, by \eqref{monotone} and Fatou's lemma, passing to the limit as $k \to\infty$ in equation \eqref{435} tells us that
$$
\int_{\{|u-\overline u|\le t\}} \big(\mathcal A(x,  \nabla u)-\mathcal A(x, \nabla \overline u)\big)\cdot (\nabla u-\nabla \overline u)\, dx=0.
$$
Thus, owing to \eqref{monotone}, we have that $\nabla u=\nabla \overline
u$ a.e. in $\{|u-\overline u|\le t\}$ for every $t>0$, whence
\begin{equation}\label{436}
\nabla u=\nabla \overline u \qquad \hbox{a.e. in }\o.
\end{equation}
%When $p\ge 2$, equation \eqref{436} immediately entails that
%$u-\overline u=c$ in $\o$ for some $c\in \R$. Indeed, since
%$u,\overline u \in V^{1, p-1}(\o)$ and $p-1\ge 1$, $u$ and
%$\overline u$ are Sobolev functions in this case.
%\par\noindent
%The case when $1<p<2$ is more delicate. Consider a family
%$\{\o_\varepsilon \}_{\varepsilon>0}$ of smooth open sets invading
%$\o$. A version of the 
By a version of the Poincar\'e inequality for domains in the class $\mathcal G_{1/\sigma '}$ \cite[Corollary 5.2.3]{Mabook} (this is also a special case of inequality \eqref{ossig}),
there exists a constant $C=C(\o)$ such that
\begin{equation}\label{437}
\left(\int_{\o }|v-{\rm med}(v)|^{\sigma '}\, dx\right)^{\frac{1}{\sigma '}}\le C 
\int_{\o}|\nabla v|\, dx,
\end{equation}
for every $v\in W^{1,1}(\o)$. Fix any $t,\tau>0$. An
application of \eqref{437} with $v= T_{\tau}(u-\T(\overline u))$,
and the use of \eqref{436} entail that
\begin{equation}\label{438}
\left( \int_{\o }|T_{\tau}(u-\T(\overline u))- {\rm med}(T_{\tau}(u-\T(\overline u)))|^{\sigma'}\, dx \right
)^{\frac{1}{\sigma'}}
\end{equation}
\begin{equation*}
\le C
\left(
\int_{\{t<|u|<t+\tau\} }|\nabla u|\, dx+
\int_{\{t-\tau<|u|<t\} }|\nabla u|\, dx
\right ).
\end{equation*}
We claim that, for each $\tau>0$, the right-hand side of
\eqref{438} converges to 0 as $t\rightarrow\infty$. To verify this
claim, choose the test function $\varphi =T_{\tau}(u_k-\T(u_k))$ in
\eqref{408bis} to deduce that
\begin{equation}\label{439}
\int_{\{t<|u_k|<t+\tau\}}B(|\nabla u_k|)\, dx\le
\int_{\{t<|u_k|<t+\tau\}}\mathcal A(x, \nabla u_k) \cdot \nabla u_k\, dx\le
\tau \int_{\{|u_k|>t\}}|f_k|\, dx.
\end{equation}
Passing to the limit as $k\rightarrow\infty$ in \eqref{439}
%,
%and making use of Fatou's lemma on the leftmost side and of the
%dominated convergence theorem on the rightmost side yield
yields, by Fatou's lemma,
\begin{equation}\label{440}
\int_{\{t<|u|<t+\tau\}}B(|\nabla u|)\, dx\le\tau
\int_{\{|u|>t\}}|f|\, dx.
\end{equation}
Hence, since the function $B$ is convex, by Jensen's inequality the first integral on the right-hand side of \eqref{438}
approaches 0 as $t\rightarrow \infty$. An analogous argument
%involving the test function $\Phi=T_{t-\tau}(u_k-\T(u_k))$
shows
that also the last integral in \eqref{438} tends to 0 as
$t\rightarrow \infty$. Since
$$
\lim_{t \to \infty} T_{\tau}(u-\T(\overline u)) -  {\rm med}(T_{\tau}(u-\T(\overline u))) =
T_{\tau}(u-\overline u) - {\rm med}(T_{\tau}(u- \overline u)) \quad \hbox{a.e. in } \o,
$$
from \eqref{438}, via Fatou's lemma, we obtain that
\begin{equation}\label{441}
\int_{\o}|T_{\tau}(u-\overline u) - {\rm med}(T_{\tau}(u- \overline u))|^{\sigma'}\,dx=0
\end{equation}
for $\tau>0$. Thus, the integrand in \eqref{441} vanishes a.e. in
$\o$ for every  $\tau >0$, and hence also  its limit
as $\tau \rightarrow \infty$ vanishes a.e. in $\o$.
Therefore, the function
%since the integrand
%in \eqref{442} approaches o a.e. in $\o$ as $\tau \rightarrow
%\infty$,
 $u-\overline u$ is constant a.e. in $\o$.

\smallskip
\par\noindent{\bf Step 7.} \emph{Equations  \eqref{umarcneu} and \eqref{gradmarcneu} hold.}
\par\noindent
Passing to the limit as $k \to\infty$ in inequality \eqref{414}, and making use of \eqref{408ter} and \eqref{425}, tell us that $u$ satisfies the assumptions of 
Lemma \ref{lemma1'}, part (ii). Hence, 
%thanks to 
%Remark \ref{rem3}, 
either inequality \eqref{4.4sig} or  inequality  \eqref{4.4infsig}  holds,
 for a suitable constant $M$ and sufficiently large $t$, depending on whether $B$ fulfills  \eqref{divinfsig} or  \eqref{convinfsig}. In the latter case, $u \in L^\infty (\o)$. In the former case,  in view of definition \eqref{Phisig},  inequality \eqref{4.4sig} yields
\begin{equation}\label{marc10}
|\{|u-{\rm med}(u)|>t\}| \Phi _\sigma(t/c)\leq c
\end{equation}
for a suitable positive constant $c$ and sufficiently large $t$.  Hence, $u \in L^{\Phi _\sigma, \infty}(\o)$. 
\\ As for $\nabla u$, if \eqref{divinfsig} is in force, then, by \eqref{4.4sig}, inequality \eqref{lem2.5} holds for a suitable $M$ and sufficicently large $s$. Therefore,
$$|\{|\nabla u|>t\}| \Psi _\sigma (t)\leq c$$
for a suitable constant $c$ and sufficiently large $t$. This implies that $\nabla u \in L^{\Psi _\sigma, \infty}(\o)$. The same conclusion follows under assumption \eqref{convinfsig}, on choosing $t=t_1$ and $s=B(t)$ in \eqref{lem2.22}, since $\lim _{t \to \infty}\phi _\sigma(t) < \infty$ in this case, and hence $\Psi_\sigma$ is equivalent to $B$  near infinity. \

\smallskip
\par\noindent{\bf Step 8.} \emph{Equation  \eqref{may1}  holds for any sequence $\{f_k\}$ as in Definition \ref{approxneu}.}
\par\noindent
This assertion follows via the arguments of Steps 2, 3 and 4, applied to the sequence $\{f_k\}$.      \qed

\medskip
\par\noindent
{\bf Proof of Theorem \ref{existneumeas}, sketched}. The proof proceeds along the same lineas as in Steps 1-5 and 7 of the proof of Theorem   \ref{existneu}. One has just to start by choosing a sequence $\{f_k\}\subset L^1(\o)\cap (W^{1,B}(\o))'$ fulfilling equation \eqref{meas6}, and such that $\int _\o f_k \, dx =0$ and 
$$\|f_k\|_{L^1(\o)} \leq 2 \|\mu \|(\o)$$
for $k \in \mathbb N$.  The function $f_k$ can be defined, for $k \in \mathbb N$, as
$f_k  = g_k - \tfrac 1{|\o|}\int _\o g_k(y)\, dy $, where $g_k : \o \to \mathbb R$ is the function given by
$$g_k(x) = \int _\o k^n \varrho (|y-x|k)\, d\mu (y) \quad \hbox{for $x \in \o$,}$$
and $\varrho : \rn \to [0, \infty]$ is a standard mollifier, namely a smooth function, compactly supported in the unit ball of $\rn$, such that $\int _\rn \varrho (x)\, dx=1$.
\\
All the equations appearing in the relevant steps then continue to hold, provided that $\|f\|_{L^1(\o)}$ is replaced by  $\|\mu \|(\o)$. 
In particular, note that, if $\varphi \in W^{1,\infty}(\o)$, then $\varphi $ is locally Lipschitz, whence  $\varphi \in C(\o)$. Moreover, $\varphi \in L^\infty (\o)$, as a consequence embedding \eqref{osnorminf}, with $B(t)=\infty$ for large values of $t$. Thus,  any function $\varphi \in W^{1,\infty}(\o)$ is admissible in \eqref{meas6}, and hence, for any such function $\varphi$, the right-hand side of \eqref{408bis} actually converges to $\int _\o \varphi \, d\mu (x)$ as $k \to \infty$.
\qed

%\bigskip
%\par\noindent
%{\bf Acknowledgment}. We wish to thank the referee for his careful reading of the manuscript and for his valuable remarks and suggestions.

\color{black}

\medskip
\par\noindent 
{\bf Acknowledgments}. 
%We wish to thank the referee for his careful reading of the manuscript and for his valuable remarks and suggestions.
%\\ \color{black}
This research was partly supported by the Research Project of the
Italian Ministry of University and Research (MIUR) Prin 2015 n.2015HY8JCC \lq\lq Partial differential equations an related analytic-geometric inequalities", by  GNAMPA   of the Italian INdAM (National Institute of High Mathematics), and by
 the Ministry of Education and Science of the Russian Federation, agreement n. 02.a03.21.0008.

%%%%%%%%%%%%%%%%%%%%%%%%%%%%%%%%%%%%%%%%%%%%%%%%%%%%%%%%%%%%%%%%%%%%%%%%%%%%%%%%%%%%%%%%%%%%%%%%%%%%%%%%
%%%%%%%%%%%%%%%%%%%%%%%%%%%%%%%%%%%%%%%%%%%%%%%%%%%%%%%%%%%%%%%%%%%%%%%%%%%%%%%%%%%%%%%%%%%%%%%%%%%%%%%%%%%

\end{document}